\documentclass[12pt,leqno]{amsart}
\usepackage{amssymb,amsmath,amsthm,bm}
\usepackage{graphicx}
\usepackage{color}
 \usepackage[textsize=tiny]{todonotes}
\usepackage[normalem]{ulem}
\usepackage[bookmarksopen,bookmarksdepth=3,colorlinks,citecolor=red,pagebackref,hypertexnames=true]{hyperref}
\usepackage[msc-links,nobysame,sorted-cites,initials]{amsrefs}
\usepackage[inline]{enumitem}
\setlist[enumerate]{font=\upshape}
\usepackage{mathtools}
\mathtoolsset{showonlyrefs}
\usepackage{tikz}
\usetikzlibrary{arrows.meta,calc}

\definecolor{darkblue}{RGB}{0,0,160}

\usepackage[lining]{libertine}
\usepackage{cabin}
\usepackage[libertine]{newtxmath}

\def\eps{\varepsilon}
\def\d{{\mathrm d}}
\def\dist{\mathrm {dist}}

\def\R {\mathbb{R}}

\def\M {{\mathrm M}}
\def\Z {{\mathbb Z}}

\def\ind{\mathbf{1}}

\newcommand{\supp}{\operatorname{supp}}

\newcommand{\pv}{\mathop{\mathrm{p.v.}}\nolimits}
\newcommand{\cW}{\mathcal W}
\newcommand{\cF}{\mathcal F}
\newcommand{\cB}{\mathcal B}
\newcommand{\cN}{\mathcal N}
\newcommand{\cD}{\mathcal D}
\newcommand{\bbS}{\mathbb S}
\newcommand{\bT}{\mathbf T}
\newcommand{\bS}{\mathbf S}
\newcommand{\sz}{\mathsf s}
\newcommand{\e}{\mathrm e}
\newcommand{\sgn}{\operatorname{sgn}}
\newcommand{\av}[2]{\langle #1\rangle_{#2}}

   \def\XXint#1#2#3{{\setbox0=\hbox{$#1{#2#3}{\int}$}
       \vcenter{\hbox{$#2#3$}}\kern-.5\wd0}}

\makeatletter
\DeclareRobustCommand\widecheck[1]{{\mathpalette\@widecheck{#1}}}
\def\@widecheck#1#2{%
    \setbox\z@\hbox{\m@th$#1#2$}%
    \setbox\tw@\hbox{\m@th$#1%
       \widehat{%
          \vrule\@width\z@\@height\ht\z@
          \vrule\@height\z@\@width\wd\z@}$}%
    \dp\tw@-\ht\z@
    \@tempdima\ht\z@ \advance\@tempdima2\ht\tw@ \divide\@tempdima\thr@@
    \setbox\tw@\hbox{%
       \raise\@tempdima\hbox{\scalebox{1}[-1]{\lower\@tempdima\box
\tw@}}}%
    {\ooalign{\box\tw@ \cr \box\z@}}}
\makeatother

\newcounter{counter}

\numberwithin{equation}{section}
\numberwithin{counter2}{section}

\newtheorem{theorem}[counter]{Theorem}

\numberwithin{withinsub}{subsection}

\numberwithin{corollary}{counter}
\newtheorem{lemma}[subsection]{Lemma}

\newtheorem{proposition}[subsection]{Proposition}

\theoremstyle{definition}

\newtheorem*{remark*}{Remark}
\newtheorem*{warn*}{A word of warning}

\newtheorem{remark}[subsection]{Remark}

\theoremstyle{definition}

 \theoremstyle{plain}

\numberwithin{figure}{section}
\numberwithin{table}{section}

\begin{document}

\title[The critical line]{The critical line for modulation invariant bilinear multipliers}

  \author[F.\ Di Plinio]{Francesco Di Plinio}
\address[F.\ Di Plinio]{Dipartimento di Matematica e Applicazioni ``R.\ Caccioppoli'' \\ \newline\indent Universit\`a di Napoli ``Federico II'', Via Cintia, Monte S.\ Angelo 80126 Napoli, Italy}
\email{\href{mailto:francesco.diplinio@unina.it}{\textnormal{francesco.diplinio@unina.it}}} \thanks{The author is partially supported by the FRA 2022--2026 Program of University of Napoli Federico II, project ReSinAPAS --
Regularity and Singularity in Analysis, PDEs, and Applied Sciences. }

\subjclass{Primary: 42B20. Secondary: 42B15}
\keywords{Bilinear Hilbert transform, modulation invariant multipliers, endpoint estimates, multi-frequency Calder\'on--Zygmund decomposition, time-frequency analysis}

\begin{abstract}
We prove that the bilinear multipliers with symbols singular along a nondegenerate line, in the class of Muscalu, Tao and Thiele, map $L^{p_1}(\mathbb R)\times L^{p_2}(\mathbb R)$ into $L^{2/3,\infty}(\mathbb R)$ on the whole critical line
\[
 \frac1{p_1}+\frac1{p_2}=\frac32,\qquad 1<p_1,p_2<2 .
\]
This is sharp, as no estimate into $L^{p,\infty}$ with $p<\frac23$ holds uniformly over the class. Moreover, the generalized restricted weak type estimates of the class extend, without logarithmic loss, to the open short sides of the hexagon of exponents. For the bilinear Hilbert transforms, this confirms a folklore conjecture recorded by Thiele and the author, and removes the doubly logarithmic term in their endpoint estimates.  The main step is an improvement over the multi-frequency Calder\'on--Zygmund decompositions already present in the literature. The bad part of this version has vanishing moments of order one at the top frequencies of the trees close to each Whitney interval, and vanishing moments of logarithmic order at the frequency bins where the tops of distant trees accumulate.
\end{abstract}
\maketitle

\section{Introduction}\label{s:intro}

The family of bilinear Hilbert transforms
\begin{equation}\label{e:bht}
 \mathrm{BHT}_b(f_1,f_2)(x)\coloneqq\pv\int_\R f_1(x-b_1t)\,f_2(x-b_2t)\,\frac{\d t}t,\qquad x\in\R,\quad b=(b_1,b_2)\in\R^2,
\end{equation}
rose to prominence from the original attempt of Calder\'on at bounding the Cauchy integral along Lipschitz curves \cites{Cal65,Cal77}, see \cite{Chr07}*{\S1.2}. The family \eqref{e:bht} displays the modulation invariance property
\begin{equation}
\label{e:lockstep}
 \mathrm{BHT}_b\big(\e^{2\pi\mathrm i\eta_1\cdot}f_1,\e^{2\pi\mathrm i\eta_2\cdot}f_2\big)=\e^{2\pi\mathrm i(\eta_1+\eta_2)\cdot}\,\mathrm{BHT}_b(f_1,f_2)\qquad\text{whenever }b_1\eta_1+b_2\eta_2=0,
\end{equation}
as the symbol $-\pi\mathrm i\sgn(b_1\xi_1+b_2\xi_2)$ of $\mathrm{BHT}_b$ is invariant under translations along the line $\{b_1\xi_1+b_2\xi_2=0\}$, across which it jumps.

The natural modulation invariant class containing the bilinear Hilbert transforms is that of Muscalu, Tao and Thiele \cite{MTT-JAMS}. Let $L\subset\R^2$ be a line through the origin which is \emph{nondegenerate}, that is, different from $\{\xi_1=0\}$, $\{\xi_2=0\}$ and $\{\xi_1+\xi_2=0\}$. The class consists of the operators
\begin{equation}\label{e:Tm}
 T_m(f_1,f_2)(x)\coloneqq\int_{\R^2}m(\xi)\,\widehat{f_1}(\xi_1)\,\widehat{f_2}(\xi_2)\,\e^{2\pi\mathrm ix(\xi_1+\xi_2)}\,\d\xi,\qquad x\in\R,
\end{equation}
whose symbols $m$ are smooth on $\R^2\setminus L$ with
\begin{equation}\label{e:symbol}
 |\partial^\gamma m(\xi)|\leq\dist(\xi,L)^{-|\gamma|},\qquad\xi\in\R^2\setminus L,\quad|\gamma|\leq\mathsf N,
\end{equation}
for a sufficiently large integer $\mathsf N$. Translating $m$ along $L$ preserves \eqref{e:symbol}, and corresponds to modulating $f_1$, $f_2$ and $T_m(f_1,f_2)$ in lockstep with \eqref{e:lockstep}. The class contains the bilinear Hilbert transforms: if $b_1,b_2\neq0$ and $b_1\neq b_2$, then
\[
 \mathrm{BHT}_b=-\pi\mathrm iT_m,\qquad m(\xi)=\sgn(b_1\xi_1+b_2\xi_2),
\]
which satisfies \eqref{e:symbol} with the nondegenerate line $L=\{b_1\xi_1+b_2\xi_2=0\}$. For the remaining parameters, $\mathrm{BHT}_b$ reduces to combinations of the Hilbert transform and of pointwise products. All constants below are allowed to depend on $L$ and $\mathsf N$.

The bilinear Hilbert transforms with nondegenerate parameters map $L^{p_1}(\R)\times L^{p_2}(\R)$ into $L^p(\R)$ in the open range
\begin{equation}\label{e:range}
 \frac1p=\frac1{p_1}+\frac1{p_2},\qquad1<p_1,p_2\leq\infty,\qquad\frac23<p<\infty,
\end{equation}
by the theorems of Lacey and Thiele \cites{LT97,LT99}, and so do the operators $T_m$, uniformly over the class, by \cite{MTT-JAMS}. While no counterexample is known to bounds for the bilinear Hilbert transforms with $p<\frac23$ and $1<p_1,p_2<\infty$, the class admits no uniform weak type bounds in that range, see \cite{Lac00}*{\S2.2}, \cite{CDPO}*{Remark~1.4} and Section~\ref{s:interp}. 

Within $1<p_1,p_2<\infty$, this article characterizes that weak type bounds hold uniformly over the class if and only if $p\geq\frac23$. 

\begin{theorem}\label{t:main}
Let $1<p_1,p_2<\infty$ and $\frac1p=\frac1{p_1}+\frac1{p_2}$. Then
\[
 \sup_m\|T_m\|_{L^{p_1}(\R)\times L^{p_2}(\R)\to L^{p,\infty}(\R)}<\infty,
\]
the supremum being over the symbols $m$ satisfying \eqref{e:symbol}, if and only if $p\geq\frac23$.
\end{theorem}

If $p>\frac23$, Theorem~\ref{t:main} is contained in \cite{MTT-JAMS}, and if $p<\frac23$, it follows from \cite{CDPO}*{Remark~1.4} by interpolation, see Section~\ref{s:interp}. The new part is the critical line $p=\frac23$, on which $1<p_1,p_2<2$. In particular, the bilinear Hilbert transforms with nondegenerate parameters map $L^{p_1}(\R)\times L^{p_2}(\R)$ into $L^{2/3,\infty}(\R)$ on the whole critical line.

The endpoint $p=\frac23$ for the bilinear Hilbert transforms was singled out in \cite{Tao07} as a model case of the logarithmic failure of time-frequency analysis. Endpoint estimates on the critical line have been formulated as generalized restricted weak type estimates for the adjoint trilinear forms, see \cite{DPT}*{(1.3)}, and the conjecture that these estimates hold on the whole critical line for the bilinear Hilbert transforms is recorded in \cite{DPT}*{Conjecture~1}. Logarithmic versions of these estimates are due to Bilyk and Grafakos \cite{BG06}. They lead to bounds in Lorentz--Orlicz spaces near the critical line \cite{CGMS09}. In \cite{DPT}*{Corollary~4}, the loss is reduced to a doubly logarithmic term, as a consequence of a mixed estimate with a general function in the first entry, \cite{DPT}*{Theorem~2}.   Proposition~\ref{p:main} below is the loss-free version of \cite{DPT}*{Theorem~2}, for the whole class, answering the question raised therein on the removal of this term, and confirming \cite{DPT}*{Conjecture~1}. The corner points of the critical line, and the quartile operator, the Walsh model of the bilinear Hilbert transform, for which the argument of Section~\ref{s:interp} proves \cite{DDP13}*{Conjecture~1}, are discussed in Remarks~\ref{r:corner} and~\ref{r:walsh}.

\subsection*{On the proof} The proof of Proposition~\ref{p:main} builds on past work of the author and collaborators \cites{DPT,DPOu}, based on the multi-frequency Calder\'on--Zygmund decomposition of Nazarov, Oberlin and Thiele \cite{NOT}; see \cite{DF} for a smooth variant. As in \cite{DPT}, the model sums representing $\Lambda_m$ are organized into forests by the size $\mathsf{s}$ of $f_3$, and $f_1\in L^p$, $1<p<2$, is split relative to the tops of the trees as \[f_1=g+\sum_Hb_H\] where $H$ ranges over the Whitney intervals of the exceptional set of $f_1$, the $L^2$-norm of the good part $g$ is controlled, for $\|f_1\|_p=|F_3|=1$, roughly as in
\begin{equation}
\label{e:expintro}
 \left\|g\right\|_2 \lesssim 
\mathsf{s}^{-2\left(\frac{1}{p}-\frac12 \right)}
\end{equation}
 and $b_H$ is supported near $H$, with vanishing moments at the frequencies relevant on $H$. The difficulty lies in the interaction of $b_H$ with wave packets localized near, but not on, $H$. In \cite{DPT}, this interaction is estimated tile by tile, with moments at about $\log N$ frequencies near each of $N$ relevant tops and with wave packets of almost exponential decay, whence the logarithmic losses. In \cite{DPOu}, moments of high order and a separation of scales trade the logarithm for a loss in the exponent \eqref{e:expintro}. Here, the interaction is summed on each Whitney interval, with two devices.
\begin{enumerate}[label=(\alph*),leftmargin=2em]
\item \emph{Near trees}, whose top interval, dilated by a bounded factor, contains $H$. At their top frequencies, $b_H$ has vanishing moments of order one. A tri-tile of relative scale $\ell(I_s)/\ell(H)\simeq2^\mu$ then interacts with $b_H$ only through a Taylor remainder of relative size $2^{-2\mu}$, while $\lesssim2^\mu$ first frequency intervals of that scale share a frequency bin of length $\ell(H)^{-1}$. The sum over the relative scales converges in each bin, and the bins are summed by orthogonality.
\item \emph{Far trees}, whose wave packets reach $H$ only through their tails, with weights which decay in the relative distance and do not depend on the number of scales of the tree. When the tops of many far trees share a bin, their tails add up to a function which is band-limited at the scale of $H$, and which is approximated on $H$ to accuracy $r^{-1}$, $r$ the total weight of the bin, by polynomials of degree $O(\log r)$ times the character of the bin. The function $b_H$ is given these moments, at the center of the bin only. As $\log r\lesssim r$, the dimension spent on the far trees is controlled by the $L^1$ norm of the counting function of the tops.
\end{enumerate}
The pairing with $f_1\ind_H$ interpolates between the orthogonality of the bins in $L^2$ and their number in $L^\infty$, and the dimension of the space of moments enters only through the Nikolskii inequality for translation invariant spaces \cite{DPOu}*{Lemma~A.2}, see also \cite{BE06}, as in the projection lemma \cite{DPOu}*{Lemma~8.1}.

\subsection*{Organization} Section~\ref{s:interp} states the main proposition, Proposition~\ref{p:main}, and deduces Theorem~\ref{t:main} from it. Section~\ref{s:model} recalls the model sums and trees of \cite{MTT-JAMS}, with the sizes of \cite{DF} applied component by component, states the level estimate, Proposition~\ref{p:level}, and deduces Proposition~\ref{p:main} from it. Section~\ref{s:lemmas} contains auxiliary lemmas, Section~\ref{s:proof} proves Proposition~\ref{p:level}, and Section~\ref{s:remarks} contains remarks.

\subsection*{Constants} The implicit constants in $\lesssim$ and $\simeq$ may depend on the parameters fixed in the statement being proved, such as the exponents, the parameter $\kappa$ of the rank one collections and the order of decay $M$ of the wave packets of Section~\ref{s:model}, but not on the functions, the collections of tri-tiles, or the intervals involved. Constants which are raised to a variable power, and thus cannot be hidden in $\lesssim$, are denoted by $C_\kappa$. Such a constant depends only on $\kappa$, and its value may change from line to line.

\section{The main proposition and the proof of Theorem~\ref{t:main}}\label{s:interp}

The operators $T_m$ are studied through their adjoint trilinear forms
\[
 \Lambda_m(f_1,f_2,f_3)\coloneqq\int_\R T_m(f_1,f_2) f_3 .
\]
The main result of this article is the following mixed estimate at the points of the critical line, with a general function in the first entry. As customary
\[\|h\|_{r,\infty}\coloneqq\sup_{\lambda>0}\lambda|\{|h|>\lambda\}|^{\frac 1r}, \qquad 0<r<\infty.\]
A major subset of a set $F$ of finite measure is a measurable set $F'\subset F$ with $|F|\leq2|F'|$.

\begin{proposition}\label{p:main}
Let $\frac12<\alpha_1<1$, $\alpha_2\coloneqq\frac32-\alpha_1$ and $p\coloneqq\frac1{\alpha_1}$. There is a constant $C$, depending only on $L$, $\mathsf N$ and $\alpha_1$, such that the following holds. For every $f_1\in L^p(\R)$ and all sets $F_2,F_3\subset\R$ of finite measure there exists a major subset $F_3'\subset F_3$, depending on $f_1,F_2,F_3$ but not on $m$, such that
\[
 |\Lambda_m(f_1,f_2,f_3)|\leq C\|f_1\|_{L^p(\R)}|F_2|^{\alpha_2}|F_3|^{-\frac12}\qquad\text{for all }|f_2|\leq\ind_{F_2},\ |f_3|\leq\ind_{F_3'},
\]
for all $m$ satisfying \eqref{e:symbol}. If $|f_1|\leq\ind_{F_1}$, then $F_3'$ may be chosen depending only on $F_1,F_2,F_3$, with $|F_1|^{\alpha_1}$ in place of $\|f_1\|_{L^p(\R)}$.
\end{proposition}

For the bilinear Hilbert transforms, Proposition~\ref{p:main} is \cite{DPT}*{Theorem~2} without the correction factor. Its case $|f_1|\leq\ind_{F_1}$ and the symmetry of the class under permutations of $f_1,f_2,f_3$ give the generalized restricted weak type estimates of \cite{DPT}*{(1.3)} on the open short sides of the hexagon of exponents, confirming \cite{DPT}*{Conjecture~1}. Proposition~\ref{p:main} is proved in Sections~\ref{s:model} to~\ref{s:proof}. The rest of this section deduces Theorem~\ref{t:main} from Proposition~\ref{p:main} at two points of the critical line, by the multilinear interpolation of Christ \cite{Chr85} and Janson \cite{Jan88}, see \cite{BOS09}*{Proposition~2.3}. The argument needs estimates without loss at two distinct points of the critical line, and it does not apply to estimates whose correction depends on the measure of the level set of the output, as in \cite{DPT}*{Theorem~2}.

 We include the short proof of the interpolation in the case needed here. It is the argument of \cite{Chr85}, with the triangle inequality replaced by the $r$-convexity of $L^{r,\infty}$, which is due to Stein, Taibleson and Weiss \cite{STW81}, see \cite{BOS09}*{(2.2)}.

\begin{lemma}\label{l:superpos}
Let $0<r<1$. For finitely many measurable functions $h_n$,
\[
 \Big\|\sum_nh_n\Big\|_{r,\infty}^r\leq\frac{2-r}{1-r}\sum_n\|h_n\|_{r,\infty}^r .
\]
\end{lemma}
\begin{proof}
Let $s_n\coloneqq\|h_n\|_{r,\infty}$ and $\lambda>0$. The set $E\coloneqq\bigcup_n\{|h_n|>\lambda\}$ has measure at most $\lambda^{-r}\sum_ns_n^r$, and off $E$ we have \[\Big|\sum_nh_n\Big|\leq\sum_n|h_n|\ind_{\{|h_n|\leq\lambda\}},\] with
\[
 \int|h_n|\ind_{\{|h_n|\leq\lambda\}}\leq\int_0^\lambda|\{|h_n|>\mu\}|\,\d\mu\leq\int_0^\lambda s_n^r\mu^{-r}\,\d\mu=\frac{s_n^r\lambda^{1-r}}{1-r} .
\]
Chebyshev's inequality on $\R\setminus E$ gives $|\{|\sum_nh_n|>\lambda\}|\leq\lambda^{-r}\frac{2-r}{1-r}\sum_ns_n^r$.
\end{proof}

\begin{lemma}\label{l:interp}
Let $0<r<1$ and $1<p_1,p_2<\infty$ with $\frac1{p_1}+\frac1{p_2}=\frac1r$, let $0<\eta\leq\min\{\frac1{p_1},\frac1{p_2}\}$, and let $T$ be a bilinear operator, defined on pairs of simple functions, such that for both choices of the sign
\begin{equation}\label{e:twopoints}
 \|T(f_1,f_2)\|_{r,\infty}\leq K|F_1|^{\frac1{p_1}\pm\eta}|F_2|^{\frac1{p_2}\mp\eta},\qquad |f_j|\leq\ind_{F_j},\ j=1,2 .
\end{equation}
Then, with $C$ depending only on $r,p_1,p_2$,
\[
 \|T(f_1,f_2)\|_{r,\infty}\leq C\,(r\eta)^{-\frac1r}K\,\|f_1\|_{p_1}\|f_2\|_{p_2}.
\]
\end{lemma}
\begin{proof}
For a simple function $f$ with decreasing rearrangement $f^*$ let
\[
 f_{(j)}\coloneqq f\,\ind_{\{f^*(2^j)<|f|\leq f^*(2^{j-1})\}},\qquad j\in\Z .
\]
These functions have disjoint supports, sum to $f$, and only finitely many of them are nonzero. Moreover $|\{f_{(j)}\neq0\}|\leq2^j$ and $\|f_{(j)}\|_\infty\leq f^*(2^{j-1})$. Decompose $f_1=\sum_jf_{1,(j)}$, $f_2=\sum_kf_{2,(k)}$, and set
\[
 A_j\coloneqq2^{j/p_1}f_1^*(2^{j-1}),\qquad B_k\coloneqq2^{k/p_2}f_2^*(2^{k-1}),
\]
so that $\|A\|_{\ell^{p_1}}\lesssim\|f_1\|_{p_1}$ and $\|B\|_{\ell^{p_2}}\lesssim\|f_2\|_{p_2}$. By homogeneity, \eqref{e:twopoints} with the upper sign if $j\leq k$ and the lower sign if $j>k$ gives
\[
 \|T(f_{1,(j)},f_{2,(k)})\|_{r,\infty}\leq K\,f_1^*(2^{j-1})2^{j(\frac1{p_1}\pm\eta)}f_2^*(2^{k-1})2^{k(\frac1{p_2}\mp\eta)}=K A_jB_k2^{-\eta|j-k|} .
\]
By Lemma~\ref{l:superpos}, Young's inequality for the convolution on $\Z$ with $w_m\coloneqq2^{-r\eta|m|}$, and H\"older's inequality with the exponents $\frac{p_1}r$ and $\frac{p_2}r$, whose reciprocals sum to $1$,
\[
 \|T(f_1,f_2)\|_{r,\infty}^r\lesssim K^r\sum_{j,k}A_j^rB_k^rw_{j-k}\leq K^r\|w\|_{\ell^1}\|A\|_{\ell^{p_1}}^r\|B\|_{\ell^{p_2}}^r .
\]
As $\|w\|_{\ell^1}\lesssim(r\eta)^{-1}$, the lemma follows.
\end{proof}

\begin{proof}[Proof of Theorem~\ref{t:main}]
If $p>\frac23$, the strong type bounds of \cite{MTT-JAMS}*{Theorem~1.1} are uniform over the class. Let $p=\frac23$, so that $1<p_1,p_2<2$, and fix $\eta>0$ such that $\frac1{p_j}\pm\eta\in(\frac12,1)$. Then \eqref{e:twopoints} holds for $T=T_m$ and $r=\frac23$. Given $|f_j|\leq\ind_{F_j}$, $\lambda>0$ and $c\in\{\pm1,\pm\mathrm i\}$, apply Proposition~\ref{p:main} with $\alpha_1=\frac1{p_1}\pm\eta$ and $F_3\coloneqq\{\mathrm{Re}(cT_m(f_1,f_2))>\lambda\}$, which has finite measure as $T_m(f_1,f_2)\in L^2(\R)$ by \cite{MTT-JAMS}*{Theorem~1.1}, and with $f_3\coloneqq c\ind_{F_3'}$, to get
\[
 \tfrac\lambda2|F_3|\leq\mathrm{Re}\,\Lambda_m(f_1,f_2,f_3)\leq C|F_1|^{\alpha_1}|F_2|^{\alpha_2}|F_3|^{-\frac12},
\]
and take the union over $c$. Lemma~\ref{l:interp} gives the bound of Theorem~\ref{t:main} for simple functions, and the extension to $L^{p_1}(\R)\times L^{p_2}(\R)$ is by density.

Finally, let $p<\frac23$, and suppose that the supremum in Theorem~\ref{t:main} is finite. Let $x_1\coloneqq(\frac1{p_1},\frac1{p_2})$, and let $x_2,x_3\in(0,1)^2$ be two points with coordinates summing to less than $1$, such that $x_1,x_2,x_3$ are not collinear. The strong type bounds of \cite{MTT-JAMS}*{Theorem~1.1} at $x_2$ and $x_3$, and the weak type bound at $x_1$, hold uniformly over the class. Hence so do the strong type bounds on the open triangle with vertices $x_1,x_2,x_3$, by the multilinear Marcinkiewicz interpolation theorem \cite{GLLZ}*{Theorem~1.1, Corollary~1.1}, which allows target exponents below $1$. This triangle contains points $(\frac1{q_1},\frac1{q_2})$ with $\frac1{q_1}+\frac1{q_2}>\frac32$, at which uniform strong type bounds are excluded by \cite{CDPO}*{Remark~1.4}.
\end{proof}

\begin{remark}[The corner $L^1\times L^2$]\label{r:corner}
The constants obtained here for Proposition~\ref{p:main}, and hence for Theorem~\ref{t:main}, blow up as the point $(\frac1{p_1},\frac1{p_2})$ approaches either corner $(1,\frac12)$ or $(\frac12,1)$ of the critical line, and we do not track their growth. At the corner $L^1\times L^2$, the known estimates carry logarithmic factors. For the bilinear Hilbert transforms, Bilyk and Grafakos \cite{BG06} prove a restricted weak type estimate at the corner up to a logarithmic factor, and \cite{DPT}*{Corollary~9} gives bounds from $L^1(\R)\times L^2(\R)$, with a logarithmic bump on the first input, into a logarithmically enlarged weak $L^{2/3}$ space. For the whole class, the sparse bounds of \cite{DF}*{Theorem~C} recover the first estimate of \cite{DPT}*{Corollary~9}, and they quantify the approach to the corner along the line $p_2=2$, from the interior of the range \eqref{e:range}:
\[
 \|T_m(f_1,f_2)\|_{L^{\frac2{3-2\varepsilon},\infty}(\R)}\lesssim\frac1\varepsilon\|f_1\|_{L^{\frac1{1-\varepsilon}}(\R)}\|f_2\|_{L^2(\R)},\qquad0<\varepsilon<2^{-6},
\]
see \cite{DF}*{Corollary~C.1}. All these estimates reach the corner with a target space which is either strictly larger than $L^{2/3,\infty}$ or of exponent strictly above $\frac23$. The improvement brought by Theorem~\ref{t:main} lies in the approach along the critical line: the target is exactly $L^{2/3,\infty}$ at every point of the open critical line, with constants which, however, do not survive at the corner.

Unfortunately, we do not know how to verify whether any of these estimates is sharp, even within the class. The example of \cite{Lac00}*{\S2.2}, realized by symbols of the class in \cite{CDPO}*{Remark~1.4}, bounds the norms of $T_m$ from below by a power $n^{\frac1{p_1}+\frac1{p_2}-\frac32}$ of a parameter $n$, which remains bounded on the whole critical line, corners included. At the corner itself, the strong type bound fails already for the bilinear Hilbert transforms. If $b_1,b_2\neq0$ and $b_1\neq b_2$, let
\[
 f_1\coloneqq\ind_{[0,1]}\in L^1(\R),\qquad f_2(y)\coloneqq\frac{\ind_{\{|y|>\e\}}(y)}{|y|^{1/2}\log|y|}\in L^2(\R).
\]
For $|x|$ large, the integrand in \eqref{e:bht} is supported on an interval of length $|b_1|^{-1}$ containing $x/b_1$, on which it does not change sign and $x-b_2t=(1-\frac{b_2}{b_1})x+O(1)$. Therefore $|\mathrm{BHT}_b(f_1,f_2)(x)|\simeq|x|^{-3/2}(\log|x|)^{-1}$ for $|x|$ large, and $\mathrm{BHT}_b(f_1,f_2)\notin L^{2/3}(\R)$. However, $\mathrm{BHT}_b(f_1,f_2)\in L^{2/3,\infty}(\R)$: its level sets $\{|\mathrm{BHT}_b(f_1,f_2)|>\lambda\}$ have measure $\lesssim\lambda^{-2/3}$ for $\lambda\leq1$ by the asymptotics above, and for $\lambda>1$ by Chebyshev's inequality, as $\mathrm{BHT}_b(f_1,f_2)\in L^2(\R)$ by \eqref{e:range} and $f_1,f_2\in L^4(\R)$. We do not know of an example ruling out the weak type bound at the corner. In particular, we do not know whether the logarithmic factors at the corner are necessary, and whether the constants of Theorem~\ref{t:main} must blow up at the corners.
\end{remark}

\begin{remark}[The quartile operator]\label{r:walsh}
For the quartile operator $V_{\mathbf{Qt}}$ and $\frac1{p_1}+\frac1{p_2}=\frac32$, $p_1,p_2>1$, \cite{DDP13}*{Theorem~2} gives
\[
 \|V_{\mathbf{Qt}}(f_1,f_2)\|_{2/3,\infty}\lesssim(p_1)'\|f_1\|_{p_1}|F_2|^{\frac1{p_2}},\qquad|f_2|\leq\ind_{F_2} .
\]
Restricted to $|f_1|\leq\ind_{F_1}$ at two points $\frac1{p_1}\pm\eta$, this is \eqref{e:twopoints}, and Lemma~\ref{l:interp} yields \cite{DDP13}*{Conjecture~1}, $V_{\mathbf{Qt}}:L^{p_1}\times L^{p_2}\to L^{2/3,\infty}$, which contains \cite{DDP13}*{Theorem~3}.
\end{remark}

\section{Model sums and the level estimate}\label{s:model}

As customary \[\langle f,g\rangle\coloneqq\int f\bar g,
\] and we normalize
\[\widehat f(\xi)\coloneqq\int f(x)\e^{-2\pi\mathrm ix\xi}\,\d x.\] The symbols $\ell(I)$ and $c(I)$ stand respectively for the length and the center of an interval $I$. Local norms are denoted as usual by \[\av{f}{p,I}\coloneqq\left( \frac{1}{|I|}\int_I|f|^p\right)^{\frac1p}, \qquad 0<p<\infty.\] The Hardy--Littlewood maximal operator is $\M$, and $\M_pf\coloneqq(\M|f|^p)^{\frac1p}$.

\subsection{Tiles and model sums}\label{ss:tiles} We use the tiles and wave packets of \cite{DF}*{\S2}. Fix two dual dyadic grids $\cD,\cD'$ and the tiling $\bbS\coloneqq\bbS_{\cD,\cD'}$ of \cite{DF}*{\S2.1}, whose elements are the \emph{tiles}
\[ P=I_P\times\omega_P\in\cD\times\cD', \qquad |I_P||\omega_P|=1.\] For $\omega\in\cD'$ and an integer $\kappa\geq0$, let $\omega^{\mathsf p(\kappa)}$ be the interval of $\cD'$ of length $2^\kappa\ell(\omega)$ containing $\omega$. Fix a large integer $M$. For a tile $P$, the class $\Phi(P)$ of \cite{DF}*{(2.1)} consists of the wave packets
\[
 \varphi(x)=\frac{\e^{2\pi\mathrm ic(\omega_P)x}}{|I_P|}\,\vartheta\Big(\frac{x-c(I_P)}{|I_P|}\Big),\qquad\supp\widehat\varphi\subset\omega_P,
\]
where
\[
 \sup_{y\in\R}\,(1+y^2)^{\frac M2}|\vartheta^{(a)}(y)|\leq1,\qquad0\leq a\leq M .
\]
In particular,
\begin{equation}\label{e:packets}
 \supp\widehat\varphi\subset\omega_P,\qquad|\varphi(x)|\lesssim|I_P|^{-1}\Big(1+\frac{\dist(x,I_P)}{|I_P|}\Big)^{-M},\qquad\varphi\in\Phi(P).
\end{equation}
The wave packet transform of \cite{DF}*{(2.2)} is \[W[f](P)\coloneqq\sup_{\varphi\in\Phi(P)}|\langle f,\varphi\rangle|.\] For an integer $\kappa\geq1$, we use the order relations on tiles of \cite{DF}*{\S7.1}:
\begin{itemize}[leftmargin=2em]
\item $P\lesssim_\kappa P'$ if $I_P\subset I_{P'}$ and $\omega_{P'}^{\mathsf p(\kappa)}\subset\omega_P^{\mathsf p(\kappa)}$;
\item $P\lesssim'_\kappa P'$ if $P\lesssim_\kappa P'$ and $P\not\lesssim_1P'$.
\end{itemize}

A \emph{tri-tile} is a triple $s=(s_1,s_2,s_3)\in\bbS^3$ of tiles with a common spatial interval $I_s$. Following \cite{MTT-JAMS}*{\S\S6--7}, in the dyadic form of \cite{DF}*{\S7.1}, a finite collection $\bS$ of tri-tiles has \emph{rank one}, with parameter an integer $\kappa\geq10$, if 
\begin{itemize}[leftmargin=2em]
\item the scales of $\bS$ are $2^{5\kappa}$-separated, that is $\ell(I_{s'})\geq2^{5\kappa}\ell(I_s)$ whenever $\ell(I_{s'})>\ell(I_s)$, $s,s'\in \bS$;
\item the maps $s\mapsto s_j$, $j=1,2,3$, are injective on $\bS$;
\item whenever $s\neq s'$ and $s_j\lesssim_1s'_j$ for some $j$,
\begin{equation}\label{e:rank1}
 s_k\lesssim_\kappa s'_k\text{ for }k=1,2,3,\qquad s_k\lesssim'_\kappa s'_k\text{ for at least two indices }k.
\end{equation}
\end{itemize}
Compare to \cite{MTT-JAMS}*{Lemma~6.2} and to the properties r1.\ to r4.\ of \cite{DF}*{\S7.1}. Subcollections of rank one collections also have rank one. Associate to a rank one collection $\bS$ the \emph{model sums}\begin{equation}\label{e:model}
 \Lambda_\bS(f_1,f_2,f_3)\coloneqq\sum_{s\in\bS}\eps_s|I_s|\prod_{j=1}^3\langle f_j,\varphi_{s_j}\rangle,\qquad|\eps_s|\leq1,\quad\varphi_{s_j}\in\Phi(s_j).
\end{equation}
For $m$ satisfying \eqref{e:symbol}, the model sum reduction of \cite{MTT-JAMS}, summarized in \cite{DF}*{\S7.1} in the present dyadic setting, writes $\Lambda_m$ as an absolutely convergent average of model sums over rank one collections, with $\kappa$ depending only on the slope of the singular line $L$, and with $M$ large if $\mathsf N$ is large. Hence, as in \cite{DPT}*{Remark~4.1}, Proposition~\ref{p:main} follows from the same estimate for $\Lambda_\bS$, uniform in $\bS$, $\eps_s$ and $\varphi_{s_j}$, provided that the major subset $F_3'$ is chosen independently of these.

\subsection{Trees and sizes} Let $\bS$ be a rank one collection of tri-tiles. For $j\in\{1,2,3\}$, a \emph{$j$-tree} is a set $\bT\subset\bS$ containing a tri-tile $s_\bT$, its top, such that $s_j\lesssim_1s_{\bT,j}$ for all $s\in\bT$, compare \cite{MTT-JAMS}*{Definition~7.1}. A \emph{tree} is a $j$-tree for some $1\leq j \leq 3$. We write
\[
 I_\bT\coloneqq I_{s_\bT},\qquad\ell_\bT\coloneqq\ell(I_\bT),\qquad\xi_{\bT,k}\coloneqq c(\omega_{s_{\bT,k}}),\qquad\bT_k\coloneqq\{s_k:s\in\bT\}\subset\bbS,
\]
the first and the last being referred to as the \emph{top interval} and the \emph{$k$-th component} of $\bT$. If $\bS'\subset\bS$ is the disjoint union of a collection $\cF$ of trees, we write 
\begin{equation}
\label{e:Nf}
\cN(\cF)\coloneqq\sum_{\bT\in\cF}|I_\bT|.
\end{equation}
We recall the trees of tiles of \cite{DF}*{\S3.1}. A set of tiles $T\subset\bbS$ is a \emph{$\kappa$-tree} with \emph{top data} $(I_T,\xi_T)\in\cD\times\R$, and \emph{top interval} $I_T$, if
\[
 I_P\subset I_T,\qquad\xi_T\in\omega_P^{\mathsf p(\kappa)}\qquad\forall P\in T,
\]
and a $\kappa$-tree $T$ is \emph{lacunary} if
\[
 \omega_P\cap\omega_{P'}=\varnothing\qquad\forall P,P'\in T\text{ with }\omega_P\neq\omega_{P'} .
\]
The next lemma shows that the components of trees of tri-tiles are $\kappa$-trees, so that the definitions and results of \cite{DF} apply component by component.

\begin{lemma}\label{l:trees}
Let $\bT$ be a $j$-tree. For $k\in\{1,2,3\}$, the component $\bT_k$ is a $\kappa$-tree with top data $(I_\bT,\xi_{\bT,k})$ in the sense of \cite{DF}*{\S3.1}, and $\bT_k$ is lacunary if $k\neq j$. Moreover:
\begin{enumerate}[label=(\roman*),leftmargin=2em]
\item at most $2^\kappa$ elements of $\bT$ share their spatial interval, and
\begin{equation}\label{e:topfreq}
 |\xi_{\bT,1}-c(\omega_{s_1})|\leq2^\kappa|\omega_{s_1}|\qquad\text{for all }s\in\bT ;
\end{equation}
\item every subset of $\bT$ is the disjoint union of $j$-trees whose top intervals are pairwise disjoint and contained in $I_\bT$.
\end{enumerate}
\end{lemma}
\begin{proof}
For $s\in\bT$, we have $I_s\subset I_\bT$ and, by \eqref{e:rank1} if $s\neq s_\bT$, $s_k\lesssim_\kappa s_{\bT,k}$, so that
\[
 \xi_{\bT,k}\in\omega_{s_{\bT,k}}^{\mathsf p(\kappa)}\subset\omega_{s_k}^{\mathsf p(\kappa)} .
\]
Thus $\bT_k$ is a $\kappa$-tree with top data $(I_\bT,\xi_{\bT,k})$. The first claim of (i) follows from the injectivity of $s\mapsto s_1$ and from \cite{DF}*{(3.1)}, by which at most $2^\kappa$ tiles of a $\kappa$-tree share their spatial interval, and \eqref{e:topfreq} follows from the case $k=1$ of the last display, as $|\omega_{s_1}^{\mathsf p(\kappa)}|=2^\kappa|\omega_{s_1}|$.
Let $k\neq j$, and suppose that $\omega_{s'_k}\subsetneq\omega_{s_k}$ for some $s,s'\in\bT$. Then $|I_{s'}|\geq2^{5\kappa}|I_s|$, so that
\[
 \xi_{\bT,k}\in\omega_{s'_k}^{\mathsf p(\kappa)}\subset\omega_{s_k} .
\]
The dyadic intervals $\omega_{s_{\bT,k}}^{\mathsf p(1)}$ and $\omega_{s_k}^{\mathsf p(1)}$ both contain $\xi_{\bT,k}$, and the former is not longer than the latter, so $s_k\lesssim_1s_{\bT,k}$. As $|I_s|<|I_{s'}|\leq|I_\bT|$, we have $s\neq s_\bT$, and $s_j\lesssim_1s_{\bT,j}$. Then \eqref{e:rank1} gives $s_i\lesssim'_\kappa s_{\bT,i}$, hence $s_i\not\lesssim_1s_{\bT,i}$, for two indices $i$, which must be the two indices different from $j$, a contradiction.
Finally, if $s,s'\in\bT$ and $I_s\subset I_{s'}$, the dyadic intervals $\omega_{s_j}^{\mathsf p(1)}$ and $\omega_{s'_j}^{\mathsf p(1)}$ contain $\omega_{s_{\bT,j}}^{\mathsf p(1)}$, and $|\omega_{s'_j}|\leq|\omega_{s_j}|$, so that $\omega_{s'_j}^{\mathsf p(1)}\subset\omega_{s_j}^{\mathsf p(1)}$ and $s_j\lesssim_1s'_j$. Hence a subset $\bT'\subset\bT$ is the disjoint union of the $j$-trees
\[
 \{s\in\bT':I_s\subset I_{s_0}\}\quad\text{with top }s_0,
\]
where $s_0\in\bT'$ is chosen once for each maximal element of $\{I_s:s\in\bT'\}$, which proves (ii).
\end{proof}

By Lemma~\ref{l:trees}, the sizes of \cite{DF}*{\S3.4} apply to the components of trees. We use the size of \cite{DF}*{(3.5)}, defined for a nonnegative function $F$ on $\bbS$ and a $\kappa$-tree $T$ by
\[
 \mathsf{size}_{2,\star,\kappa}(F,T)\coloneqq\sup_{U} \Big(\frac1{|I_U|}\sum_{P\in U}|I_P|\,F(P)^2\Big)^{\frac12},
\]
the supremum being taken over the lacunary $\kappa$-trees $U\subset T$, with top interval $I_U$. For $\bS'\subset\bS$, $k\in\{1,2,3\}$ and $f\in L^1(\R)+L^\infty(\R)$, let
\[
 \sz_k(f;\bS')\coloneqq\sup\big\{\mathsf{size}_{2,\star,\kappa}(W[f],\bT_k):\ \bT\subset\bS'\text{ is a $j$-tree for some }j\neq k\big\}.
\]
As $\{s\}$ is a $j$-tree for every $j$, and as $\bT_k$ is a lacunary $\kappa$-tree with top interval $I_\bT$ if $\bT$ is a $j$-tree with $j\neq k$,
\begin{equation}\label{e:single}
 W[f](s_k)\leq\sz_k(f;\bS'),\qquad\sum_{s'\in\bT}|I_{s'}|\,W[f](s'_k)^2\leq|I_\bT|\,\sz_k(f;\bS')^2,
\end{equation}
for $s\in\bS'$ and for $j$-trees $\bT\subset\bS'$ with $j\neq k$.
The bound of \cite{DF}*{Proposition~3.11}, applied with $J=I_\bT$ to the components $\bT_k$, gives
\begin{equation}\label{e:sizeM}
 \sz_k(f;\bS')\lesssim\sup_{s\in\bS'}\inf_{x\in I_s}\M f(x).
\end{equation}
If $\bT\subset\bS'$ is a $j$-tree, bounding the $j$-th factor by the first estimate of \eqref{e:single}, and the other two by the Cauchy--Schwarz inequality and the second bound, yields the single tree estimate, compare \cite{DF}*{(7.2)},
\begin{equation}\label{e:tree}
 \sum_{s\in\bT}|I_s|\prod_{k=1}^3W[f_k](s_k)\leq|I_\bT|\prod_{k=1}^3\sz_k(f_k;\bS') .
\end{equation}
Similarly, choosing $k\in\{2,3\}\setminus\{j\}$ and bounding by \eqref{e:single} the factor of the other index in $\{2,3\}$ gives
\[
 \sum_{s\in\bT}|I_s|\,W[f_2](s_2)^2\,W[f_3](s_3)^2\leq|I_\bT|\,\sz_2(f_2;\bS')^2\,\sz_3(f_3;\bS')^2 .
\]
Hence, if $\bS'=\bigsqcup_{\bT\in\cF}\bT$ is a disjoint union of trees,
\begin{equation}\label{e:l2coef}
 \sum_{s\in\bS'}|I_s|\,W[f_2](s_2)^2\,W[f_3](s_3)^2\leq\cN(\cF)\,\sz_2(f_2;\bS')^2\,\sz_3(f_3;\bS')^2 .
\end{equation}

The forests are produced by the size lemma for rank one collections, \cite{MTT-JAMS}*{Lemma~7.7}, see also \cite{DPT}*{Lemma~5.2}; its single tiling counterpart is the embedding \cite{DF}*{Proposition~3.14}. These are stated for sizes defined through trees and one wave packet per tile, and they apply to $\sz_k$: by Lemma~\ref{l:trees}(ii), the supremum over lacunary subtrees $U\subset\bT_k$ in the definition of $\sz_k$ may be restricted to $U=\bT_k$, and, $\bS$ being finite, $W[f]$ may be replaced by $|\langle f,\varphi_P\rangle|$ for wave packets $\varphi_P\in\Phi(P)$ attaining it up to a factor $2$.

\begin{lemma}\label{l:size}
Let $k\in\{1,2,3\}$, $f\in L^2(\R)$, $\bS'\subset\bS$ and $\sigma\geq\sz_k(f;\bS')$. Then $\bS'=\bS^{\mathrm{lo}}\sqcup\bS^{\mathrm{hi}}$, where $\sz_k(f;\bS^{\mathrm{lo}})\leq\frac\sigma2$ and $\bS^{\mathrm{hi}}$ is the disjoint union of a collection $\cF$ of trees with $\cN(\cF)\lesssim\sigma^{-2}\|f\|_2^2$.
\end{lemma}

The trees of \cite{MTT-JAMS} need not contain their top, and those of \cite{DPT} are defined through top data; by the argument for Lemma~\ref{l:trees}(ii), each of them is a disjoint union of trees as above, with pairwise disjoint top intervals contained in its own. Combining Lemma~\ref{l:size} with the single tree estimate \eqref{e:tree} yields the following forest estimate, compare \cite{DPT}*{Lemma~5.3}.

\begin{proposition}\label{p:forest}
Let $A>0$, $f_1\in L^2(\R)$ and $f_2,f_3\in L^1(\R)+L^\infty(\R)$. If $\bS'\subset\bS$ is the disjoint union of a collection $\cF$ of trees with
\[
 \cN(\cF)\leq A^2\,\sz_3(f_3;\bS')^{-2},
\]
then
\[
 \sum_{s\in\bS'}|I_s|\prod_{k=1}^3W[f_k](s_k)\lesssim A\,\|f_1\|_2\,\sz_2(f_2;\bS') .
\]
\end{proposition}
\begin{proof}
Let $\sigma_k\coloneqq\sz_k(f_k;\bS')$. We may assume that $\|f_1\|_2=1$ and $\sigma_1,\sigma_3>0$, as otherwise the left hand side vanishes by \eqref{e:single}. Let $n_0\in\Z$ be the least integer with $A\sigma_1\leq2^{n_0}\sigma_3$. If $n_0\leq0$, then \eqref{e:tree} bounds the left hand side by
\[
 \cN(\cF)\,\sigma_1\sigma_2\sigma_3\leq A^2\frac{\sigma_1}{\sigma_3}\,\sigma_2\leq A\sigma_2 .
\]
If $n_0>0$, iterating Lemma~\ref{l:size} with $f=f_1$ for $n=0,\dots,n_0-1$ produces forests $\bS'_n$, and the remaining collection $\bS'_{n_0}$ is the disjoint union of the collection $\cF_{n_0}$ of trees given by Lemma~\ref{l:trees}(ii) for the sets $\bT\cap\bS'_{n_0}$, $\bT\in\cF$, so that $\cN(\cF_{n_0})\leq\cN(\cF)\leq A^2\sigma_3^{-2}\leq2^{2n_0}\sigma_1^{-2}$. This decomposes
\[
 \bS'=\bigsqcup_{n=0}^{n_0}\bS'_n,\qquad\bS'_n=\bigsqcup_{\bT\in\cF_n}\bT,\qquad\cN(\cF_n)\lesssim\frac{2^{2n}}{\sigma_1^2},\qquad\sz_1(f_1;\bS'_n)\leq\frac{\sigma_1}{2^n} .
\]
Applying \eqref{e:tree} to the trees of each $\cF_n$, and as $2^{n_0-1}\sigma_3<A\sigma_1$ by the choice of $n_0$, the left hand side is at most
\[
 \sum_{n=0}^{n_0}\cN(\cF_n)\,\frac{\sigma_1}{2^n}\,\sigma_2\sigma_3\lesssim\sum_{n=0}^{n_0}2^n\frac{\sigma_3}{\sigma_1}\,\sigma_2\lesssim2^{n_0}\frac{\sigma_3}{\sigma_1}\,\sigma_2\lesssim A\sigma_2 .\qedhere
\]
\end{proof}

\subsection{The level estimate} Proposition~\ref{p:main} rests on the following estimate for a single forest, proved in Section~\ref{s:proof}.

\begin{proposition}\label{p:level}
Let $1<p_1<2$, $\alpha\coloneqq\frac{2}{p_1}-1$ and $C_0>0$. Let $f_1\in L^{p_1}(\R)$ with $\|f_1\|_{p_1}\leq1$. Let $\bS=\bigsqcup_{\bT\in\cF}\bT$ be a rank one collection which is the disjoint union of a collection $\cF$ of trees, such that
\[I_s\not\subset E_1\coloneqq\left\{\M_{p_1}f_1>C_0\right\} \qquad \forall s\in\bS.\] Then
\[
 |\Lambda_\bS(f_1,f_2,f_3)|\lesssim\max\left\{1,  \cN(\cF) \right\}^{\frac{\alpha+1}{2}}\sz_2(f_2;\bS)\,\sz_3(f_3;\bS)
\]
for all $f_2,f_3\in L^2(\R)$, with a constant depending only on $p_1$, $C_0$, $\kappa$ and $M$.
\end{proposition}

For comparison, as $|f_1|\leq C_0$ outside $E_1$, Proposition~\ref{p:forest} bounds the contribution of $f_1\ind_{\R\setminus E_1}$ by $\cN(\cF)^{1/2}\sz_2\sz_3$, while a tree by tree estimate of the rest of $f_1$ gives $\cN(\cF)\sz_2\sz_3$. The point of Proposition~\ref{p:level} is the exponent $\frac{1+\alpha}2$, with no loss in $\cN(\cF)$.

\subsection{Proof of Proposition~\ref{p:main}}\label{ss:proofmixed} Recall that $p=\frac1{\alpha_1}$, and let $\alpha\coloneqq\frac2p-1\in(0,1)$, so that
\[
 \alpha_1=\frac{1+\alpha}2,\qquad\alpha_2=1-\frac\alpha2,\qquad\alpha_2-\frac1{p'}=\frac12 .
\]
We may assume that $|F_2|>0$ and $\|f_1\|_p=1$. If $|F_2|\geq|F_3|$, let $F_3'\coloneqq F_3$; the strong type bound \cite{MTT-JAMS}*{Theorem~1.1} at the exponents $(\frac1p,\frac1{2p'},\frac1{2p'})$ gives
\[
 |\Lambda_m(f_1,f_2,f_3)|\lesssim|F_2|^{\frac1{2p'}}|F_3|^{\frac1{2p'}}\leq|F_2|^{\frac1{p'}}\leq|F_2|^{\alpha_2}|F_3|^{-\frac12} .
\]
Let $|F_2|<|F_3|$. By Section~\ref{ss:tiles}, it suffices to bound $\Lambda_\bS$ for a rank one collection $\bS$, with $F_3'$ independent of $\bS$, and by the dyadic scaling invariance of the model sums, compare \cite{DPT}*{Remark~4.2}, we may assume that $|F_3|\simeq1$, with implicit constants depending only on $\kappa$. For a large constant $C$, let
\[
 E_1\coloneqq\{\M_pf_1>C\},\qquad E\coloneqq E_1\cup\left\{\M\frac{\ind_{F_2}}{|F_2|}>C\right\},\qquad F_3'\coloneqq F_3\setminus\bar E,
\]
where $\bar E\coloneqq\{\M\ind_E>\frac15\}$. The set $\bar E$ contains $3Q$ for every interval $Q\subset E$, so that $F_3'$ does not depend on the dyadic grids. As $|F_3|\simeq1$, these are, up to the values of the constants, the sets of \cite{DPT}*{(6.1)--(6.2)} for the tuple $(\alpha_1,\alpha_2,-\frac12)$, $h_1=f_1$ and $h_2=\ind_{F_2}$, and $F_3'$ is a major subset of $F_3$ if $C$ is large, as
\[
 |\bar E|\leq10|E|\lesssim C^{-p}+C^{-1}.
\]

Let $\bS^{\mathrm{in}}\coloneqq\{s\in\bS:I_s\subset E\}$ and $\bS^{\mathrm{out}}\coloneqq\bS\setminus\bS^{\mathrm{in}}$. The proof of \cite{DPT}*{Proposition~6.1} uses the model sums only through \eqref{e:packets} and the injectivity of $s\mapsto s_j$, and gives
\[
 |\Lambda_{\bS^{\mathrm{in}}}(f_1,f_2,f_3)|\lesssim|F_2|^{\alpha_2},\qquad|f_2|\leq\ind_{F_2},\ |f_3|\leq\ind_{F_3'} .
\]
As $I_s\not\subset E$ for $s\in\bS^{\mathrm{out}}$, \eqref{e:sizeM} gives
\[
 \sz_2(f_2;\bS^{\mathrm{out}})\lesssim|F_2|,\qquad\sz_3(f_3;\bS^{\mathrm{out}})\lesssim1 .
\]
Let $n_0\geq0$ be the least integer with $2^{-2n_0}\leq|F_2|$, so that $2^{-n_0}\simeq|F_2|^{1/2}$ as $|F_2|\lesssim1$. Iterating Lemma~\ref{l:size} with $f=f_3$ decomposes $\bS^{\mathrm{out}}$, up to tri-tiles with $W[f_3](s_3)=0$, as
\[
 \bS^{\mathrm{out}}=\bigsqcup_{k\geq0}\bS_k,\qquad\bS_k=\bigsqcup_{\bT\in\cF_k}\bT,\qquad\cN(\cF_k)\lesssim2^{2k},\qquad\sz_3(f_3;\bS_k)\lesssim2^{-k},
\]
see \cite{DPT}*{(6.3)--(6.5)}. For $k>n_0$, as $\sz_2(f_2;\bS_k)\lesssim|F_2|\lesssim2^{-n_0}|F_2|^{1/2}$, iterating Lemma~\ref{l:size} with $f=f_2$ for $n=n_0,\dots,k-1$, and organizing the leftover tri-tiles into trees by Lemma~\ref{l:trees}(ii), decomposes further
\[
 \bS_k=\bigsqcup_{n=n_0}^k\bS_{k,n},\qquad\bS_{k,n}=\bigsqcup_{\bT\in\cF_{k,n}}\bT,\qquad\cN(\cF_{k,n})\lesssim2^{2n},\qquad\sz_2(f_2;\bS_{k,n})\lesssim2^{-n}|F_2|^{\frac12},
\]
see \cite{DPT}*{(6.8)}. As $E_1\subset E$, Proposition~\ref{p:level} with $p_1=p$ and $C_0=C$ gives
\[
 \sum_{k\leq n_0}|\Lambda_{\bS_k}(f_1,f_2,f_3)|\lesssim\sum_{k\leq n_0}2^{(1+\alpha)k}|F_2|\,2^{-k}\lesssim\frac{2^{\alpha n_0}}\alpha|F_2|\simeq\frac{|F_2|^{1-\frac\alpha2}}\alpha
\]
and, as $\sz_3(f_3;\bS_{k,n})\lesssim2^{-k}$,
\[
 \sum_{k>n_0}\sum_{n=n_0}^k|\Lambda_{\bS_{k,n}}(f_1,f_2,f_3)|\lesssim\sum_{k>n_0}\sum_{n=n_0}^k2^{\alpha n}|F_2|^{\frac12}2^{-k}\lesssim\frac{2^{-(1-\alpha)n_0}}{\alpha(1-\alpha)}|F_2|^{\frac12}\simeq\frac{|F_2|^{1-\frac\alpha2}}{\alpha(1-\alpha)} .
\]
As $\alpha_2=1-\frac\alpha2$ and $|F_3|\simeq1$, this proves the first statement. The second follows from the same proof, applied to $|F_1|^{-\alpha_1}f_1$ with $|F_1|^{-\alpha_1}\ind_{F_1}$ in place of $f_1$ in the definition of $E_1$, which makes $F_3'$ depend only on $F_1,F_2,F_3$. \qed

\section{Four lemmas}\label{s:lemmas}

Throughout, $1<p<2$ and $\alpha=\frac2p-1\in(0,1)$. For an interval $J$, let
\[
 \chi_J^\infty(x)\coloneqq\left(1+\frac{|x-c(J)|}{\ell(J)}\right)^{-1},\qquad x\in\R ,
\]
and, for $K\geq0$, define the tailed and the localized local norms
\[
 \av{f}{\infty,J,+K}\coloneqq\big\|(\chi_J^\infty)^{K}f\big\|_\infty,\qquad\av{f}{\infty,J,-K}\coloneqq\big\|(\chi_J^\infty)^{-K}f\big\|_\infty ;
\]
the latter is finite only if $f$ decays at order $K$ away from $J$.
 First, recall the projection lemma of \cite{DPOu}*{Lemma~8.1}, stated for a general translation invariant space as allowed by its proof, \cite{DPOu}*{Lemma~A.1}, and in terms of the local norms of Section~\ref{s:model}.

\begin{lemma}\label{l:proj}
Let $H$ be an interval and $V$ a linear space of continuous functions on $\R$ of finite linear dimension $d$, invariant under translations. Let $\{e_i\}$ be an orthonormal basis of $\{v\ind_{3H}:v\in V\}$ in $L^2(3H)$, and let
\[
 \Pi u\coloneqq\ind_{3H}\sum_i\langle u,e_i\rangle e_i,\qquad u\in L^1(3H),
\]
be the orthogonal projection onto this space. Let $u\in L^p(\R)$ vanish outside $H$. Then $b\coloneqq u-\Pi u$ vanishes outside $3H$, $\langle b,v\rangle=0$ for all $v\in V$, and
\begin{align}
 \label{e:projpair}\frac{|\langle u,\Pi f\rangle|}{|H|}&\lesssim d^{\frac\alpha2} \av{u}{p,H}\av{f}{2,3H},\qquad f\in L^2(3H),\\
 \label{e:projgood}\av{\Pi u}{2,3H}&\lesssim d^{\frac\alpha2}\av{u}{p,H} .
\end{align}
\end{lemma}
\begin{proof}
The Nikolskii inequality for translation invariant spaces \cite{DPOu}*{Lemma~A.2}, see also \cite{BE06}, gives $\|v\|_{L^\infty(H)}\lesssim d^{\frac12}\av{v}{2,3H}$ for $v\in V$. As $1-\frac2{p'}=\alpha$, H\"older's inequality on $H$ yields
\[
 \av{v}{p',H}\leq\|v\|_{L^\infty(H)}^{\alpha}\,\av{v}{2,H}^{1-\alpha}\lesssim d^{\frac\alpha2}\av{v}{2,3H},\qquad v\in V .
\]
As $\Pi f$ coincides on $H$ with an element of $V$, and $\av{\Pi f}{2,3H}\leq\av{f}{2,3H}$,
\[
 |\langle u,\Pi f\rangle|\leq|H|\,\av{u}{p,H}\,\av{\Pi f}{p',H}\lesssim d^{\frac\alpha2}\,|H|\,\av{u}{p,H}\,\av{f}{2,3H},
\]
which is \eqref{e:projpair}. Then \eqref{e:projgood} is \eqref{e:projpair} with $f=\Pi u$, as \[|3H|\av{\Pi u}{2,3H}^2=\langle u,\Pi u\rangle.\] Finally $\langle b,e_i\rangle=0$ for all $i$.
\end{proof}

\begin{lemma}\label{l:bern}
Let $K,q\geq0$ be integers, $J$ an interval and $b\geq\ell(J)^{-1}$. If $F\in L^1(\R)$ and $\supp\widehat F\subset[-b,b]$, then
\[
 \av{F^{(q)}}{\infty,J,+K}\lesssim(8\pi b)^q\av{F}{\infty,J,+K},
\]
with implicit constant depending only on $K$.
\end{lemma}
\begin{proof}
Pick $\phi$ smooth, $\phi=1$ on $[-1,1]$ and $\supp\phi\subset[-2,2]$. 
Then \[
F^{(q)}=F*\theta_b^{(q)}, \qquad \widehat{\theta_b} \coloneqq \phi \left({\textstyle\frac {\cdot}{b}}\right).\] Integrating by parts $K+2$ times in the Fourier integral defining $\theta_b^{(q)}$ gives \[|\theta_b^{(q)}(y)|\lesssim(1+q)^{K+2}(4\pi b)^q\,b\,(1+b|y|)^{-K-2}\lesssim(8\pi b)^q\,b\,(1+b|y|)^{-K-2}.\]  
As $b\ell(J)\geq1$, for all $x,y\in\R$
\[
\chi_J^\infty(x-y)^{-1}\leq1+\frac{|x-c(J)|}{\ell(J)}+\frac{|y|}{\ell(J)}\leq\chi_J^\infty(x)^{-1}(1+b|y|),
\]
whence $|F(x-y)|\leq\av{F}{\infty,J,+K}\,\chi_J^\infty(x)^{-K}(1+b|y|)^K$. Therefore
\[
\chi_J^\infty(x)^K|F^{(q)}(x)|\leq\av{F}{\infty,J,+K}\int(1+b|y|)^K|\theta_b^{(q)}(y)|\,\d y\lesssim(8\pi b)^q\av{F}{\infty,J,+K},
\]
as claimed.
\end{proof}
Next, denote by $\cB_H\coloneqq\{\ell(H)^{-1}[k,k+1): k\in \mathbb Z\}$ the \emph{frequency bins} associated to the interval $H$.

\begin{lemma}\label{l:bins}
Let $H$ be an interval and $\cB\subset\cB_H$ be a finite subcollection. For each $B\in\cB$, let there be given a function $g_B$ on $4H$ which is smooth at the scale of $H$, in the sense that
\[
 \left\|g_B^{(q)}\right\|_{L^\infty(4H)}\leq A_B\,\ell(H)^{-q},\qquad q=0,1,2,
\]
for some $A_B>0$. Then
\[
 \left\langle\sum_{B\in\cB}\e^{2\pi\mathrm ic(B)\cdot}g_B\right\rangle_{2,3H}\lesssim\left(\sum_{B\in\cB}A_B^2\right)^{\frac12} .
\]
\end{lemma}
\begin{proof}
Fix a smooth $\eta_H$ with $\ind_{3H}\leq\eta_H\leq\ind_{4H}$ and $|\eta_H^{(q)}|\lesssim\ell(H)^{-q}$, and let $G_B\coloneqq\eta_H\e^{2\pi\mathrm ic(B)\cdot}g_B$ for $B\in\cB$. For $B\neq B'$, integrating by parts twice, 
\[
\begin{split}
 \left|\langle G_B,G_{B'}\rangle\right|&=\left|\int\e^{2\pi\mathrm i(c(B)-c(B'))x}\,\eta_H(x)^2\,g_B(x)\,\overline{g_{B'}(x)}\,\d x \right|  \lesssim\frac{|H|\,A_BA_{B'}}{\big(1+\ell(H)|c(B)-c(B')|\big)^2} .
 \end{split}
\]
As $\ell(H)|c(B)-c(B')|$ ranges over the integers, Schur's test gives
\[
 \left\langle\sum_{B\in\cB}\e^{2\pi\mathrm ic(B)\cdot}g_B\right\rangle_{2,3H}^2\leq \frac{1}{|3H|}\Big\|\sum_{B\in\cB}G_B\Big\|_2^2\lesssim \sum_{B\in\cB}A_B^2 ,
\]
which is the claimed bound.
\end{proof}

The last lemma  interpolates between orthogonality in $L^2$ and a count of frequencies in $L^\infty$. Let $O\subset\R$ be an open set of finite measure, and $\cW$ a partition of $O$ into intervals, such as its Whitney intervals. A function $\Upsilon\geq0$ on $O$, constant on each $H\in\cW$, plays the role of a counting function, and a function $F$ with $|F|\leq\Upsilon$ that of a sum of at most $\Upsilon$ terms of unit size. 

\begin{lemma}\label{l:sum}
Let $O$ and $\cW$ be as above, $u\in L^p(O)$   and $\Upsilon\geq0$ constant on each $H\in\cW$. If $|F|\leq\Upsilon$ on $O$, then
\[
 \|uF\|_{L^1(O)}\leq|O|^{\frac{1-\alpha}2}\,\|\Upsilon\|_{L^1(O)}^{\alpha}\,\|F\|_{L^2(O)}^{1-\alpha} \left( \sup_{H \in \mathcal W}  \av{u}{p,H} \right).
\]
\end{lemma}
\begin{proof} By scaling $u$, we can assume  $\av{u}{p,H}\leq1$ for each $H \in \mathcal W$.
As $|F|\leq\Upsilon$ and $\frac1p+\frac{1-\alpha}2=1$, H\"older's inequality gives
\[
 \|uF\|_{L^1(O)}\leq\big\|u\,\Upsilon^\alpha|F|^{1-\alpha}\big\|_{L^1(O)}\leq\|u\Upsilon^\alpha\|_{L^p(O)}\,\|F\|_{L^2(O)}^{1-\alpha} .
\]
Write $\Upsilon_H$ for the value of $\Upsilon$ on $H\in\cW$. As $\av{u}{p,H}\leq1$ and $\alpha p<1$, H\"older's inequality once more gives
\[
 \|u\Upsilon^\alpha\|_{L^p(O)}^p=\sum_{H\in\cW}\Upsilon_H^{\alpha p}\int_H|u|^p\leq\int_O\Upsilon^{\alpha p}\leq|O|^{1-\alpha p}\,\|\Upsilon\|_{L^1(O)}^{\alpha p},
\]
and the lemma follows, as $\frac1p-\alpha=\frac{1-\alpha}2$.
\end{proof}
\section{Proof of Proposition~\ref{p:level}}\label{s:proof}

Let $p_1,\alpha,C_0,f_1,E_1$ and $\bS=\bigsqcup_{\bT\in\cF}\bT$  be as in the statement of Proposition~\ref{p:level}, and write $f\coloneqq f_1$. The lemmas of Section~\ref{s:lemmas} will be applied with $p\coloneqq p_1$. By scaling  $f_2,f_3\in L^2(\R)$, we can reduce to the case 
\begin{equation}
\label{e:norm37}
\max\left\{\sz_2(f_2;\bS),\sz_3(f_3;\bS)\right\} \leq 1.
\end{equation}
It is convenient to use the representation
\begin{equation}
\label{e:rep37}
\Lambda_\bS(h,f_2,f_3) = \langle h,\Psi\rangle, \qquad \Psi\coloneqq\sum_{s\in\bS}\overline{\beta_s}\,\varphi_{s_1}, \qquad 
 \beta_s\coloneqq\eps_s|I_s|\langle f_2,\varphi_{s_2}\rangle\langle f_3,\varphi_{s_3}\rangle.
\end{equation} Write $\cN\coloneqq\cN(\cF)$, see \eqref{e:Nf}.
 As $|\langle f_j,\varphi_{s_j}\rangle|\leq W[f_j](s_j)$ for $j\in\{2,3\}$, the normalization \eqref{e:norm37}, the bounds \eqref{e:single}, \eqref{e:l2coef} and Proposition~\ref{p:forest} with $A=\cN^{\frac12}$ give
\begin{equation}\label{e:coef}
 |\beta_s|\leq|I_s|,\qquad\sum_{s\in\bS}\frac{|\beta_s|^2}{|I_s|}\leq\cN,\qquad\|\Psi\|_2\lesssim\cN^{\frac12} .
\end{equation}
The wave packets $\varphi_{s_1}$, $s\in\bS$, in \eqref{e:rep37} have frequency supports in the intervals $\omega_{s_1}$, and are localized to $I_s$: by \eqref{e:packets},
\begin{equation}\label{e:chi}
 \av{\varphi_{s_1}}{\infty,I_s,-M}\lesssim|I_s|^{-1},\qquad s\in\bS .
\end{equation}
In particular, the finite sum $\Psi$ is integrable and bounded. By Lemma~\ref{l:trees}, at most $2^\kappa$ tri-tiles of a tree share their spatial interval.

\subsection{Whitney intervals, near and far trees}\label{ss:setup} The set $E_1$ is open and $|E_1|\lesssim1$ by the maximal theorem. Let $\cW$ be the collection of the maximal dyadic intervals $H$ with $9H\subset E_1$. These intervals are pairwise disjoint with union $E_1$, the intervals $3H$ have bounded overlap, and
\begin{equation}\label{e:whitney}
 \dist(H,\R\setminus E_1)\geq4\ell(H),\qquad\av{f}{p,H}\lesssim1 .
\end{equation}
For  each $H\in\cW$,  split
\[
\begin{split}
&\mathbf{S}= \mathbf{S}^{\mathsf{coarse}}(H) \sqcup  \mathbf{S}^{\mathsf{fine}}(H),\\ & \mathbf{S}^{\mathsf{coarse}}(H) \coloneqq \left\{ s\in \mathbf S:  \ell(I_s)\geq\ell(H) \right\}, \quad  
\mathbf{S}^{\mathsf{fine}}(H) \coloneqq \left\{ s\in \mathbf S:  \ell(I_s)<\ell(H) \right\}.
\end{split}
\]
Since $I_s\not\subset E_1$, \eqref{e:whitney} gives
\begin{equation}\label{e:finesep}
 \dist(H,I_s)\geq3\ell(H)\qquad\text{for every }s\in\bS^{\mathsf{fine}}(H) .
\end{equation}
For each $H\in\cW$, split also the trees into \emph{near}, \emph{far} and \emph{fine} trees,
\[
\begin{split}
&\cF=\cF^{\mathsf{near}}(H)\sqcup\cF^{\mathsf{far}}(H)\sqcup\cF^{\mathsf{fine}}(H),\\
&\cF^{\mathsf{near}}(H)\coloneqq\{\bT\in\cF:\ \ell_\bT\geq\ell(H),\ \dist(H,I_\bT)\leq8\ell_\bT\},\\
&\cF^{\mathsf{far}}(H)\coloneqq\{\bT\in\cF:\ \ell_\bT\geq\ell(H),\ \dist(H,I_\bT)>8\ell_\bT\}, \\
 &\cF^{\mathsf{fine}}(H)\coloneqq\{\bT\in\cF:\ \ell_\bT<\ell(H)\}.
\end{split}
\]
Every tree of $\cF^{\mathsf{fine}}(H)$ is contained in $\bS^{\mathsf{fine}}(H)$. The near trees are counted by
\[
 n_H\coloneqq\#\cF^{\mathsf{near}}(H).
\]
A far tree $\bT\in\cF^{\mathsf{far}}(H)$ reaches $H$ only through the tails of its wave packets. Summing \eqref{e:chi}, with \eqref{e:coef}, over the tri-tiles of $\bT$, first at each scale and then over the scales, bounds the contribution of $\bT\cap\bS^{\mathsf{coarse}}(H)$ to $\Psi$ on $4H$ by a constant times the \emph{weight}
\[
 w_\bT(H)\coloneqq\Big(\frac{\ell_\bT}{\dist(H,I_\bT)}\Big)^{M-1}\ll 1,
\]
see the proof of Lemma~\ref{l:far}. Accordingly, the far trees are counted with these weights, sorted by the frequency bin $B\in\cB_H$ of Lemma~\ref{l:bins} containing their top frequency:
\[
 n^{\mathsf{far}}_{H,B}\coloneqq\sum_{\substack{\bT\in\cF^{\mathsf{far}}(H)\\ \xi_{\bT,1}\in B}}w_\bT(H),\qquad n^{\mathsf{far}}_H\coloneqq\sum_{\bT\in\cF^{\mathsf{far}}(H)}w_\bT(H)=\sum_{B\in\cB_H}n^{\mathsf{far}}_{H,B} .
\]
Thus $n^{\mathsf{far}}_{H,B}$ bounds the contribution on $4H$ of the far trees with top frequency in $B$, see \eqref{e:farder}, and it determines the number of moments of $b_H$ at the bin $B$ in Section~\ref{ss:decomp}. The total, weighted, count of the trees relevant to $H$ is
\[
 \Upsilon_H\coloneqq n_H+n^{\mathsf{far}}_H ,
\]
and its integral over $E_1$ is controlled by $\cN$. 
\begin{lemma}\label{l:count} With the previous definitions, 
\begin{equation}\label{e:count}
 \sum_{H\in\cW}|H|\,\Upsilon_H\lesssim\cN .
\end{equation}
\end{lemma}
\begin{proof}
Write the left hand side of \eqref{e:count} as a sum over the trees,
\[
 \sum_{H\in\cW}|H|\,\Upsilon_H=\sum_{\bT\in\cF}\Bigg(\sum_{\substack{H\in\cW\\ \bT\in\cF^{\mathsf{near}}(H)}}|H|\;+\sum_{\substack{H\in\cW\\ \bT\in\cF^{\mathsf{far}}(H)}}|H|\,w_\bT(H)\Bigg),
\]
so that it suffices to bound each of the two inner sums by a constant times $|I_\bT|$. If $\bT\in\cF^{\mathsf{near}}(H)$, then $\ell(H)\leq\ell_\bT$ and $\dist(H,I_\bT)\leq8\ell_\bT$, so that $H\subset19I_\bT$. As the intervals of $\cW$ are pairwise disjoint, the first inner sum is at most $19|I_\bT|$. If $\bT\in\cF^{\mathsf{far}}(H)$, then $\ell(H)\leq\ell_\bT<\dist(H,I_\bT)$, so that
\[
 8\ell_\bT<\dist(H,I_\bT)\leq\dist(x,I_\bT)\leq\dist(H,I_\bT)+\ell(H)<2\dist(H,I_\bT),\qquad x\in H,
\]
and therefore
\[
 |H|\,w_\bT(H)\leq2^{M-1}\int_H\Big(\frac{\ell_\bT}{\dist(x,I_\bT)}\Big)^{M-1}\d x .
\]
The intervals $H$ with $\bT\in\cF^{\mathsf{far}}(H)$ are pairwise disjoint and contained in $\{x:\dist(x,I_\bT)>8\ell_\bT\}$, so that the second inner sum is at most
\[
 2^{M-1}\int_{\dist(x,I_\bT)>8\ell_\bT}\Big(\frac{\ell_\bT}{\dist(x,I_\bT)}\Big)^{M-1}\d x\lesssim\ell_\bT=|I_\bT| .
\]
Summing over $\bT\in\cF$ proves \eqref{e:count}.
\end{proof}
\subsection{The decomposition}\label{ss:decomp} For $H\in\cW$, the function $f\ind_H$ is split into a good part $g_H$, its projection as in Lemma~\ref{l:proj} onto a translation invariant space $V_H$ of finite dimension $d_H$, and a bad part $b_H$, which vanishes outside $3H$ and has vanishing moments, in the sense that $\langle b_H,v\rangle=0$ for all $v\in V_H$. The space $V_H$ is dictated by the coarse tri-tiles, sorted by the frequency bin $B\in\cB_H$ containing the top frequency of their tree. On $4H$, the moments must remove from the contribution of each bin to $\Psi$ a function of $V_H$, leaving a remainder which is smooth at the scale of $H$ in the sense of Lemma~\ref{l:bins}, with $A_B\lesssim1$, so that the bins can be summed by that lemma. This is achieved in two steps, for the near and for the far trees.

\subsubsection*{Near trees} Let $s$ be a coarse tri-tile of a tree $\bT\in\cF^{\mathsf{near}}(H)$, with $\ell(I_s)\simeq2^\mu\ell(H)$. On $4H$, the wave packet $\varphi_{s_1}$, demodulated by the top frequency $\xi_{\bT,1}$, is smooth at the scale $\ell(I_s)$. Hence it differs from its Taylor polynomial of order one at $c(H)$ by a remainder of relative size $2^{-2\mu}$, and the remainders are summable over the relative scales $2^\mu$ within each bin, see Lemma~\ref{l:near}. The Taylor polynomials, times the characters of the top frequencies, belong to
\[
 V^{\mathsf{near}}_H\coloneqq\operatorname{span}\big\{x^{z}\e^{2\pi\mathrm i\xi_{\bT,1}x}:\ \bT\in\cF^{\mathsf{near}}(H),\ {z}\in\{0,1\}\big\} .
\]

\subsubsection*{Far trees} Let $B\in\cB_H$. On $4H$, the coarse tri-tiles of the far trees with top frequency in $B$ contribute to $\Psi$ the character of $c(B)$ times a function which is band-limited at the scale of $H$ and of size $n^{\mathsf{far}}_{H,B}$, see Lemma~\ref{l:far}. If $n^{\mathsf{far}}_{H,B}\gg 1$, this size is too large for Lemma~\ref{l:bins}, but the function differs from its Taylor polynomial at $c(H)$ of degree $\sim 1+\log n^{\mathsf{far}}_{H,B}$ by a bounded remainder. These Taylor polynomials, times the characters of the centers of the bins, belong to
\[
\begin{split}
 V^{\mathsf{far}}_H&\coloneqq\operatorname{span}\big\{x^q\e^{2\pi\mathrm ic(B)x}:\ B\in\cB_H,\ n^{\mathsf{far}}_{H,B}>1,\ 0\leq q\leq Q_{H,B}\big\},\\
 Q_{H,B}&\coloneqq\big\lceil A(1+\log n^{\mathsf{far}}_{H,B})\big\rceil,
\end{split}
\]
where $A$ is a large constant depending only on $\kappa$.

Let $V_H\coloneqq V^{\mathsf{near}}_H+V^{\mathsf{far}}_H$. For all $\xi\in\R$ and $Q\geq0$, the span of the functions $x^q\e^{2\pi\mathrm i\xi x}$, $0\leq q\leq Q$, is invariant under translations, hence so is $V_H$. Let $\Pi_H$ be the projection of Lemma~\ref{l:proj} for $H$ and $V=V_H$, and
\[
 g_H\coloneqq\Pi_H(f\ind_H),\qquad b_H\coloneqq f\ind_H-g_H,\qquad G\coloneqq\sum_{H\in\cW}g_H .
\]
By Lemma~\ref{l:proj}, $g_H$ and $b_H$ vanish outside $3H$, and $b_H$ is orthogonal to $V_H$.

\begin{lemma}\label{l:good}
For $H\in\cW$, the dimension $d_H$ of $V_H$ satisfies
\begin{equation}\label{e:dim}
 d_H\lesssim\Upsilon_H .
\end{equation}
Moreover, the series defining $G$ converges in $L^2(\R)$, and
\begin{equation}\label{e:good}
 \sum_{H\in\cW}\|g_H\|_{L^2(3H)}^2\lesssim\cN^\alpha,\qquad\|G\|_2\lesssim\cN^{\frac\alpha2} .
\end{equation}
\end{lemma}
\begin{proof}
As $Q_{H,B}+1\lesssim n^{\mathsf{far}}_{H,B}$ whenever $n^{\mathsf{far}}_{H,B}>1$,
\[
 d_H\leq2n_H+\sum_{\substack{B\in\cB_H\\ n^{\mathsf{far}}_{H,B}>1}}(Q_{H,B}+1)\lesssim n_H+\sum_{B\in\cB_H}n^{\mathsf{far}}_{H,B}=\Upsilon_H .
\]
By \eqref{e:projgood}, \eqref{e:whitney} and \eqref{e:dim},
\[
 \|g_H\|_{L^2(3H)}^2=|3H|\,\av{g_H}{2,3H}^2\lesssim d_H^\alpha\,|H|\lesssim\Upsilon_H^\alpha\,|H|,\qquad H\in\cW .
\]
The intervals of $\cW$ are a pairwise disjoint cover of $E_1$. H\"older's inequality, \eqref{e:count} and $|E_1|\lesssim1$ thus give
\[
 \sum_{H\in\cW}\|g_H\|_{L^2(3H)}^2\lesssim\sum_{H\in\cW}|H|\,\Upsilon_H^\alpha\leq\Big(\sum_{H\in\cW}|H|\,\Upsilon_H\Big)^\alpha|E_1|^{1-\alpha}\lesssim\cN^\alpha .
\]
Finally, as the intervals $3H$ have bounded overlap,
\[
 \Big\|\sum_{H\in\cW}g_H\Big\|_2^2\lesssim\sum_{H\in\cW}\|g_H\|_{L^2(3H)}^2 ,
\]
which gives convergence of the series and the bound for $\|G\|_2$.
\end{proof}
Since $|f|\leq C_0$ almost everywhere outside $E_1$, we have $\|f\ind_{\R\setminus E_1}\|_2\lesssim1$. As
\[
 f=f\ind_{\R\setminus E_1}+\sum_{H\in\cW}f\ind_H=f\ind_{\R\setminus E_1}+G+\sum_{H\in\cW}b_H,
\]
the representation \eqref{e:rep37} gives
\begin{equation}\label{e:split}
 \Lambda_\bS(f,f_2,f_3)=\langle f\ind_{\R\setminus E_1},\Psi\rangle+\langle G,\Psi\rangle+\sum_{H\in\cW}\langle b_H,\Psi\rangle .
\end{equation}
Here the last series converges absolutely: $\Psi$ is integrable and bounded, so that, by H\"older's inequality and the bounded overlap of the intervals $3H$,
\[
 \sum_{H\in\cW}|\langle f\ind_H,\Psi\rangle|\leq\|f\|_p\|\Psi\|_{p'},\qquad\sum_{H\in\cW}|\langle g_H,\Psi\rangle|\lesssim\Big(\sum_{H\in\cW}\|g_H\|_{L^2(3H)}^2\Big)^{\frac12}\|\Psi\|_2 .
\]
By the Cauchy--Schwarz inequality, \eqref{e:coef} and Lemma~\ref{l:good}, the first two terms of \eqref{e:split} are bounded by
\begin{equation}\label{e:goodterms}
 \big(\|f\ind_{\R\setminus E_1}\|_2+\|G\|_2\big)\|\Psi\|_2\lesssim\big(1+\cN^{\frac\alpha2}\big)\cN^{\frac12}\lesssim\max\{1,\cN\}^{\frac{1+\alpha}2} .
\end{equation}
Therefore Proposition~\ref{p:level} follows from the next lemma, whose proof occupies the rest of the section.

\begin{lemma}\label{l:bad}
The bad parts satisfy
\begin{equation}\label{e:bad}
 \sum_{H\in\cW}|\langle b_H,\Psi\rangle|\lesssim\cN^{\frac{1+\alpha}2} .
\end{equation}
\end{lemma}

\subsection{Splitting of the bad part}\label{ss:split} For $s\in\bS$, let $\bT(s)$ be the tree of $\cF$ containing $s$. Fix $H\in\cW$. If $s\in\bS^{\mathsf{coarse}}(H)$, then $\ell_{\bT(s)}\geq\ell(I_s)\geq\ell(H)$, so that $\bT(s)$ is a near or a far tree. Accordingly, for $B\in\cB_H$ let
\begin{align*}
 \bS^{\mathrm n}_{H,B}&\coloneqq\{s\in\bS^{\mathsf{coarse}}(H):\ \bT(s)\in\cF^{\mathsf{near}}(H),\ \xi_{\bT(s),1}\in B\},\\
 \bS^{\mathrm f}_{H,B}&\coloneqq\{s\in\bS^{\mathsf{coarse}}(H):\ \bT(s)\in\cF^{\mathsf{far}}(H),\ \xi_{\bT(s),1}\in B\},
\end{align*}
and
\[
 \Psi^{\mathrm n}_{H,B}\coloneqq\sum_{s\in\bS^{\mathrm n}_{H,B}}\overline{\beta_s}\,\varphi_{s_1},\qquad\Psi^{\mathrm f}_{H,B}\coloneqq\sum_{s\in\bS^{\mathrm f}_{H,B}}\overline{\beta_s}\,\varphi_{s_1} .
\]
As $\bS^{\mathsf{coarse}}(H)$ is the disjoint union of the collections $\bS^{\mathrm n}_{H,B}$ and $\bS^{\mathrm f}_{H,B}$, $B\in\cB_H$, the coarse and the fine part of $\Psi$ are
\begin{equation}\label{e:Psisplit}
\begin{split}
 \Psi^{\mathsf{coarse}}_H&\coloneqq\sum_{s\in\bS^{\mathsf{coarse}}(H)}\overline{\beta_s}\,\varphi_{s_1}=\sum_{B\in\cB_H}\big(\Psi^{\mathrm n}_{H,B}+\Psi^{\mathrm f}_{H,B}\big),\\
 \Psi^{\mathsf{fine}}_H&\coloneqq\sum_{s\in\bS^{\mathsf{fine}}(H)}\overline{\beta_s}\,\varphi_{s_1}=\Psi-\Psi^{\mathsf{coarse}}_H .
\end{split}
\end{equation}
Hence Lemma~\ref{l:bad} follows from the two estimates
\begin{align}
 \label{e:badcoarse}\sum_{H\in\cW}|\langle b_H,\Psi^{\mathsf{coarse}}_H\rangle|&\lesssim\cN^{\frac{1+\alpha}2},\\
 \label{e:badfine}\sum_{H\in\cW}|\langle b_H,\Psi^{\mathsf{fine}}_H\rangle|&\lesssim\cN^{\frac{1+\alpha}2},
\end{align}
proved in Subsections~\ref{ss:coarse} and~\ref{ss:fine}, after the analysis of the single bins in Lemmata~\ref{l:near} and~\ref{l:far}.

Let $B\in\cB_H$ and $s\in\bS^{\mathrm n}_{H,B}\cup\bS^{\mathrm f}_{H,B}$. By definition, the top frequency $\xi_{\bT(s),1}$ belongs to $B$, an interval of length $\ell(H)^{-1}$ centered at $c(B)$. By \eqref{e:topfreq}, the frequency support $\omega_{s_1}$ of $\varphi_{s_1}$, an interval of length $\ell(I_s)^{-1}$, has center within distance $2^\kappa\ell(I_s)^{-1}$ of $\xi_{\bT(s),1}$. Hence
\begin{equation}\label{e:topbin}
 |\xi_{\bT(s),1}-c(B)|\leq\ell(H)^{-1},\qquad\omega_{s_1}\subset\xi_{\bT(s),1}+\big[-2^{\kappa+1}\ell(I_s)^{-1},2^{\kappa+1}\ell(I_s)^{-1}\big] .
\end{equation}
The following lemma is where the moments of order one at the tops enter.

\begin{lemma}\label{l:near} Let $H\in\cW$ and $B\in\cB_H$ with $\bS^{\mathrm n}_{H,B}\neq\varnothing$. There is a function $F^{\mathrm n}_{H,B}$ on $4H$ such that
\begin{equation}
\label{e:theid}
\left\langle b_H,\Psi^{\mathrm n}_{H,B}\right\rangle=\left\langle \e^{-2\pi\mathrm ic(B)\cdot} b_H,F^{\mathrm n}_{H,B}\right\rangle,\qquad\left|\left(F^{\mathrm n}_{H,B}\right)^{(z)}(x) \right|\lesssim\ell(H)^{-{z}},\qquad x\in4H,\ {z}\leq2 .
\end{equation}
\end{lemma}
\begin{proof}
Let $s\in\bS^{\mathrm n}_{H,B}$, and let
\[
 \phi_s(x)\coloneqq\e^{-2\pi\mathrm i\xi_{\bT(s),1}x}\varphi_{s_1}(x),\qquad x\in\R,
\]
be the wave packet $\varphi_{s_1}$ demodulated by the top frequency of its tree. Then $\widehat{\phi_s} $ is supported in $\omega_{s_1}-\xi_{\bT(s),1}$, hence in $[-b,b]$ with $b\coloneqq2^{\kappa+1}\ell(I_s)^{-1}$, by \eqref{e:topbin}. Let $J$ be the interval of length $4\ell(I_s)$ centered at $c(H)$. This interval is chosen so that, by the coarseness of  $s$, $4H\subset J$ and $\chi_J^\infty\geq\frac23$ on $4H$. Furthermore, by \eqref{e:chi} and the triangle inequality,
\[
\av{\phi_s}{\infty,J,+M}\lesssim\chi_{I_s}^\infty(c(H))^M\av{\varphi_{s_1}}{\infty,I_s,-M}\lesssim\frac{\chi_{I_s}^\infty(c(H))^M}{|I_s|}.
\]
As $b\geq\ell(J)^{-1}$, Lemma~\ref{l:bern} with $K=M$ applies to $\phi_s$, yielding \[
 |\phi_s^{(z)}(x)|\lesssim\frac{\chi_{I_s}^\infty(c(H))^M}{\ell(I_s)^z|I_s|},\qquad x\in4H,\ z\in \{0,1,2\}.
\]
Let $P_s$ be the Taylor polynomial of $\phi_s$ of order one at $c(H)$, and $\Delta_s\coloneqq\phi_s-P_s$ its remainder. Then Taylor's formula and the case $z=2$ of the last display give, for ${z}\in\{0,1,2\}$,
\begin{equation}\label{e:Rs}
 \left\|\Delta_s^{(z)}\right\|_{L^\infty(4H)}\leq(2\ell(H))^{2-{z}}\sup_{4H}|\phi_s''|\lesssim\ell(H)^{-{z}}\Big(\frac{\ell(H)}{\ell(I_s)}\Big)^2\frac{\chi_{I_s}^\infty(c(H))^M}{|I_s|} .
\end{equation}
We prove the identity in \eqref{e:theid}. As $\bT(s)\in\cF^{\mathsf{near}}(H)$, the function $\e^{2\pi\mathrm i\xi_{\bT(s),1}x}P_s(x)$ belongs to $V^{\mathsf{near}}_H$, is thus orthogonal to $b_H$ and may be subtracted off in the last step of the first line of  
\[
\begin{split}
 \langle b_H,\Psi^{\mathrm n}_{H,B}\rangle&=\sum_{s\in\bS^{\mathrm n}_{H,B}}\beta_s\,\big\langle b_H,\e^{2\pi\mathrm i\xi_{\bT(s),1}\cdot}\phi_s\big\rangle=\sum_{s\in\bS^{\mathrm n}_{H,B}}\beta_s\,\big\langle b_H,\e^{2\pi\mathrm i\xi_{\bT(s),1}\cdot}\left[\phi_s-P_s\right]\big\rangle\\
 &=\Big\langle\e^{-2\pi\mathrm ic(B)\cdot}b_H,\sum_{s\in\bS^{\mathrm n}_{H,B}}\overline{\beta_s}\,\e^{2\pi\mathrm i(\xi_{\bT(s),1}-c(B))\cdot}\Delta_s\Big\rangle .
\end{split}
\]
As $b_H$ vanishes outside $3H$, the identity in \eqref{e:theid} holds with
\[
 F^{\mathrm n}_{H,B}(x)\coloneqq\sum_{s\in\bS^{\mathrm n}_{H,B}}\overline{\beta_s}\,\e^{2\pi\mathrm i(\xi_{\bT(s),1}-c(B))x}\Delta_s(x),\qquad x\in4H .
\]
We now bound the derivatives of $F^{\mathrm n}_{H,B}$. By \eqref{e:topbin}, the derivative of order $k$ of $\e^{2\pi\mathrm i(\xi_{\bT(s),1}-c(B))x}$ is bounded by $(2\pi\ell(H)^{-1})^k$. Hence the Leibniz rule, \eqref{e:Rs} and $|\beta_s|\leq|I_s|$ give
\[
 \big|\big(F^{\mathrm n}_{H,B}\big)^{(z)}(x)\big|\lesssim\ell(H)^{-{z}}\sum_{s\in\bS^{\mathrm n}_{H,B}}\Big(\frac{\ell(H)}{\ell(I_s)}\Big)^2\chi_{I_s}^\infty(c(H))^M,\qquad x\in4H,\ {z}\in\{0,1,2\} .
\]
Sort the tri-tiles of $\bS^{\mathrm n}_{H,B}$ by their relative scale,
\[
 \bS_\mu\coloneqq\big\{s\in\bS^{\mathrm n}_{H,B}:\ 2^\mu\ell(H)\leq\ell(I_s)<2^{\mu+1}\ell(H)\big\},\qquad\mu\geq0 .
\]
As $\ell(H)/\ell(I_s)\sim 2^{-\mu}$ for $s\in\bS_\mu$, the claimed bound follows from
\begin{equation}\label{e:nearcount}
 \sum_{s\in\bS_\mu}\chi_{I_s}^\infty(c(H))^M\lesssim2^\mu,\qquad\mu\geq0 .
\end{equation}
As the lengths of dyadic intervals are powers of $2$, the tri-tiles of $\bS_\mu$ have the common length $\ell_\mu\coloneqq2^\mu\ell(H)$. For each $I\in\cD$ of length $\ell_\mu$, the tri-tiles $s\in\bS_\mu$ with $I_s=I$ have distinct first frequency intervals, by the injectivity of $s\mapsto s_1$, of length $\ell_\mu^{-1}=2^{-\mu}\ell(H)^{-1}$, and contained in an interval of length $\lesssim\ell(H)^{-1}$, by \eqref{e:topbin}; hence there are $\lesssim2^\mu$ of them. As the intervals $I\in\cD$ of length $\ell_\mu$ are pairwise disjoint,
\[
 \sum_{s\in\bS_\mu}\chi_{I_s}^\infty(c(H))^M\lesssim2^\mu\sum_{\substack{I\in\cD\\ \ell(I)=\ell_\mu}}\chi_I^\infty(c(H))^M\lesssim2^\mu . \qedhere
\]
\end{proof}

The tops of many far trees may lie in a single bin. The following lemma is where the moments of logarithmic order at the merged bins enter.

\begin{lemma}\label{l:far}
Let $H\in\cW$ and $B\in\cB_H$. Then
\begin{equation}\label{e:farder}
 \left| \left[ \e^{-2\pi\mathrm ic(B)\cdot}\Psi^{\mathrm f}_{H,B} \right]^{(q)}(x)\right|\lesssim\Big(\frac{C_\kappa}{\ell(H)}\Big)^q\,n^{\mathsf{far}}_{H,B},\qquad x\in4H,\ q\geq0,
\end{equation}
with implicit constant independent of $q$. Moreover, there is a function $F^{\mathrm f}_{H,B}$ on $4H$ such that
\[
 \left\langle b_H,\Psi^{\mathrm f}_{H,B}\right\rangle=\left\langle\e^{-2\pi\mathrm ic(B)\cdot}b_H,F^{\mathrm f}_{H,B}\right\rangle,\qquad\left|\left(F^{\mathrm f}_{H,B}\right)^{(z)}(x)\right|\lesssim\frac{\min\{1,n^{\mathsf{far}}_{H,B}\}}{\ell(H)^{z}},\qquad x\in4H,\ {z}\leq2 .
\]
\end{lemma}
\begin{proof}
Let
\[
 \phi_B(x)\coloneqq\e^{-2\pi\mathrm ic(B)x}\Psi^{\mathrm f}_{H,B}(x),\qquad x\in\R,
\]
be the far part of the bin $B$ demodulated by the center of $B$, so that \eqref{e:farder} is a bound for the derivatives of $\phi_B$, and $|\phi_B|=|\Psi^{\mathrm f}_{H,B}|$. As $\ell(I_s)\geq\ell(H)$ for $s\in\bS^{\mathrm f}_{H,B}$, \eqref{e:topbin} gives
\[
 \supp\widehat{\phi_B}\subset[-b,b],\qquad b\coloneqq2^{\kappa+2}\ell(H)^{-1} .
\]
As $\chi_{4H}^\infty\geq\frac23$ on $4H$, Lemma~\ref{l:bern} with $J=4H$ and $K=M-1$ reduces \eqref{e:farder} to
\begin{equation}\label{e:farweight}
 \av{\Psi^{\mathrm f}_{H,B}}{\infty,4H,+(M-1)}\lesssim n^{\mathsf{far}}_{H,B} .
\end{equation}
By the definition of $n^{\mathsf{far}}_{H,B}$, it suffices to show that the contribution
\[
 \Psi_\bT\coloneqq\sum_{s\in\bT\cap\bS^{\mathsf{coarse}}(H)}\overline{\beta_s}\,\varphi_{s_1}
\]
of a single tree $\bT\in\cF^{\mathsf{far}}(H)$ satisfies
\begin{equation}\label{e:farT}
 \av{\Psi_\bT}{\infty,4H,+(M-1)}\lesssim w_\bT(H)=\Big(\frac{\ell_\bT}{\dist(H,I_\bT)}\Big)^{M-1} ,
\end{equation}
that is, $|\Psi_\bT(y)|\,\chi_{4H}^\infty(y)^{M-1}\lesssim w_\bT(H)$ for $y\in\R$.
Recall that $\dist(H,I_\bT)>8\ell_\bT$. We begin with the contribution of a single scale. At each scale $\ell\leq\ell_\bT$ of $\bT$, the intervals $I_s$ are subintervals of $I_\bT$ of length $\ell$, each shared by at most $2^\kappa$ tri-tiles of $\bT$. Hence, by \eqref{e:chi} and \eqref{e:coef},
\begin{equation}\label{e:onescale}
 \sum_{\substack{s\in\bT\\ \ell(I_s)=\ell}}|\beta_s\varphi_{s_1}(y)|\lesssim\min\Big\{1,\Big(\frac\ell{\dist(y,I_\bT)}\Big)^{M-1}\Big\},\qquad y\in\R .
\end{equation}
We distinguish two cases. Near $H$, that is, if
$
 |y-c(H)|\leq\frac{\dist(H,I_\bT)}{2},
$
then $\dist(y,I_\bT) >4\ell_\bT$. Summing the second bound of \eqref{e:onescale} over the scales $\ell\leq\ell_\bT$, which is a geometric series, we obtain
\[
 |\Psi_\bT(y)|\,\chi_{4H}^\infty(y)^{M-1}\leq|\Psi_\bT(y)|\lesssim\Big(\frac{\ell_\bT}{\dist(y,I_\bT)}\Big)^{M-1}\lesssim\Big(\frac{\ell_\bT}{\dist(H,I_\bT)}\Big)^{M-1}=w_\bT(H) .
\]
Away from $H$, that is, if
$
 |y-c(H)|>\frac{\dist(H,I_\bT)}{2},
$
we use that  $\bT\cap\bS^{\mathsf{coarse}}(H)$ has at most $1+\log_2(\ell_\bT/\ell(H))\leq2\ell_\bT/\ell(H)$ scales, each contributing $\lesssim1$ by the first bound of \eqref{e:onescale}. Therefore, as $\ell_\bT\geq\ell(H)$,
\[
 |\Psi_\bT(y)|\,\chi_{4H}^\infty(y)^{M-1}\lesssim\frac{\ell_\bT}{\ell(H)}\Big(\frac{\ell(H)}{\dist(H,I_\bT)}\Big)^{M-1}\leq\Big(\frac{\ell_\bT}{\dist(H,I_\bT)}\Big)^{M-1}=w_\bT(H) .
\]
This proves \eqref{e:farT}, hence \eqref{e:farweight} and \eqref{e:farder}.

We prove the second claim, distinguishing two cases. Undoing the demodulation, $\Psi^{\mathrm f}_{H,B}(x)=\e^{2\pi\mathrm ic(B)x}\phi_B(x)$, so that
\begin{equation}\label{e:farid}
 \langle b_H,\Psi^{\mathrm f}_{H,B}\rangle=\big\langle\e^{-2\pi\mathrm ic(B)\cdot}b_H,\phi_B\big\rangle .
\end{equation}
If $n^{\mathsf{far}}_{H,B}\leq1$, the second claim thus holds with $F^{\mathrm f}_{H,B}\coloneqq\phi_B$, by \eqref{e:farder} with $q\leq2$.

If $n^{\mathsf{far}}_{H,B}>1$, recall from Section~\ref{ss:decomp} that $b_H$ has vanishing moments at the bin $B$ up to the degree
\[
 Q\coloneqq Q_{H,B}=\big\lceil A(1+\log n^{\mathsf{far}}_{H,B})\big\rceil,
\]
where $A$ is a large constant depending only on $\kappa$, chosen below. That is, $b_H$ is orthogonal to $\e^{2\pi\mathrm ic(B)x}P(x)$ for every polynomial $P$ of degree at most $Q$. Let $P_B$ be the Taylor polynomial of $\phi_B$ of degree $Q$ at $c(H)$, and $\Delta_B\coloneqq\phi_B-P_B$ its remainder. Subtracting $P_B$ in \eqref{e:farid} gives
\[
 \langle b_H,\Psi^{\mathrm f}_{H,B}\rangle=\big\langle\e^{-2\pi\mathrm ic(B)\cdot}b_H,\phi_B-P_B\big\rangle=\big\langle\e^{-2\pi\mathrm ic(B)\cdot}b_H,\Delta_B\big\rangle ,
\]
and, as $b_H$ vanishes outside $3H$, the identity of the second claim holds with $F^{\mathrm f}_{H,B}\coloneqq\Delta_B$ on $4H$. It remains to bound the derivatives of $\Delta_B$. As $Q\geq A\geq2$, Taylor's formula on $4H$ and \eqref{e:farder} with $q=Q+1$ give, for ${z}\in\{0,1,2\}$,
\begin{equation}
\label{e:thels}
 \big\|\Delta_B^{(z)}\big\|_{L^\infty(4H)}\leq\frac{(2\ell(H))^{Q+1-{z}}}{(Q+1-{z})!}\sup_{4H}\big|\phi_B^{(Q+1)}\big|\lesssim\ell(H)^{-{z}}\,\frac{C_\kappa^Q}{Q!}\,n^{\mathsf{far}}_{H,B} .
\end{equation}
Here factors polynomial in $Q$ have been absorbed in $C_\kappa^Q$. Now, if $A\geq\e^2C_\kappa$, then $Q\geq\e^2C_\kappa$ as well, so that $\e\,C_\kappa/Q\leq\e^{-1}$. Hence, by $k!\geq(k/\e)^k$ and as $Q\geq\log n^{\mathsf{far}}_{H,B}$,
\[
 \frac{C_\kappa^Q}{Q!}\,n^{\mathsf{far}}_{H,B}\leq\Big(\frac{\e\,C_\kappa}{Q}\Big)^Q n^{\mathsf{far}}_{H,B}\leq\e^{-Q}\,n^{\mathsf{far}}_{H,B}\leq1,
\]
and the left hand side of \eqref{e:thels} is $\lesssim\ell(H)^{-{z}}$, as claimed.
\end{proof}

\subsection{Conclusion for the coarse tri-tiles}\label{ss:coarse} We prove \eqref{e:badcoarse}. Begin by  applying Lemmata~\ref{l:near} and~\ref{l:far} for each $H\in \mathcal W$ and $B \in \mathcal B_H$, setting $F^{\mathrm n}_{H,B}\coloneqq0$ if $\bS^{\mathrm n}_{H,B}=\varnothing$. This yields
\begin{equation}
\label{e:eachHB}
 \left\|\left[F^{\mathrm n}_{H,B}+F^{\mathrm f}_{H,B}\right]^{(z)}\right\|_{L^\infty(4H)}\lesssim\ell(H)^{-z}\left(\ind_{\{\bS^{\mathrm n}_{H,B}\neq\varnothing\}}+\min\left\{1,n^{\mathsf{far}}_{H,B}\right\}\right),\qquad z\in\{0,1,2\} .
\end{equation}
As $\bS$ is finite, the right hand side vanishes for all but finitely many $B\in\cB_H$. For $H\in\cW$, let
\[
 F_H(x)\coloneqq\sum_{B\in\cB_H}\e^{2\pi\mathrm ic(B)x}\big(F^{\mathrm n}_{H,B}(x)+F^{\mathrm f}_{H,B}(x)\big),\qquad x\in4H.
\]
By virtue of
 \eqref{e:eachHB}, we may apply Lemma~\ref{l:bins} with $g_B=F^{\mathrm n}_{H,B}+F^{\mathrm f}_{H,B}$, and with $A_B$ a constant times the quantity in parentheses in \eqref{e:eachHB}. This quantity is at most $2$, at most $n_H$ bins contain the top frequency of a near tree, and $\min\{1,n^{\mathsf{far}}_{H,B}\}\leq n^{\mathsf{far}}_{H,B}$, so that
\[
 \sum_{B\in\cB_H}A_B^2\lesssim\sum_{B\in\cB_H}A_B\lesssim n_H+n^{\mathsf{far}}_H=\Upsilon_H .
\]
Therefore the triangle inequality and Lemma~\ref{l:bins} give
\begin{equation}\label{e:FH}
 \|F_H\|_{L^\infty(4H)}\lesssim\Upsilon_H,\qquad\av{F_H}{2,3H}^2\lesssim\Upsilon_H .
\end{equation}

 By \eqref{e:Psisplit}, the support conditions on $b_H$, Lemmata~\ref{l:near} and~\ref{l:far}, and the self-adjointness of $\Pi_H$,
\begin{equation}\label{e:coarseid}
 \langle b_H,\Psi^{\mathsf{coarse}}_H\rangle= \langle b_H,F_H \rangle = \langle b_H,F_H \ind_{3H}\rangle =\langle f\ind_H,F_H\rangle-\langle f\ind_H,\Pi_H(F_H\ind_{3H})\rangle .
\end{equation}
Thus \eqref{e:badcoarse} follows by summing over $H\in\cW$ the absolute values of the two terms on the right hand side of \eqref{e:coarseid}, which we bound in turn. For the first term, we apply Lemma~\ref{l:sum} on $O=E_1$ to $u=f$ and to the functions
\begin{equation}
\label{e:Fdef}
 F\coloneqq\sum_{H\in\cW}F_H\ind_H,\qquad\Upsilon\coloneqq\sum_{H\in\cW}\Upsilon_H\ind_H ,
\end{equation}
which coincide with $F_H$ and $\Upsilon_H$ on each $H\in\cW$. By the first bound in \eqref{e:FH}, $|F|\lesssim\Upsilon$ on $E_1$. By the second bound in \eqref{e:FH} and \eqref{e:count},
\[
 \|F\|_{L^2(E_1)}^2+\|\Upsilon\|_{L^1(E_1)}\lesssim\sum_{H\in\cW}|H|\,\Upsilon_H\lesssim\cN .
\]
As $|E_1|\lesssim1$, Lemma~\ref{l:sum} and \eqref{e:whitney} give
\[
 \sum_{H\in\cW}|\langle f\ind_H,F_H\rangle|\leq\|fF\|_{L^1(E_1)}\lesssim\cN^{\frac{1+\alpha}2} .
\]
For the second term, \eqref{e:projpair}, \eqref{e:whitney}, \eqref{e:dim}, \eqref{e:FH}, H\"older's inequality and \eqref{e:count} give
\[
 \sum_{H\in\cW}|\langle f\ind_H,\Pi_H(F_H\ind_{3H})\rangle|\lesssim\sum_{H\in\cW}d_H^{\frac\alpha2}\,|H|\,\av{F_H}{2,3H}\lesssim\sum_{H\in\cW}|H|\,\Upsilon_H^{\frac{1+\alpha}2}\lesssim\cN^{\frac{1+\alpha}2} .
\]
This proves \eqref{e:badcoarse}.

\subsection{Fine tri-tiles}\label{ss:fine} We prove \eqref{e:badfine}. The moments of $b_H$ play no role here. The spatial interval of a fine tri-tile is shorter than $H$ and, by \eqref{e:finesep}, lies at distance at least $3\ell(H)$ from $H$, so that its wave packet reaches $3H$ only through its tail. Measured in units of the length of the spatial interval, this distance is large because of the separation in space, and larger still because of the separation in scale. Accordingly, fix $H\in\cW$, and sort the fine tri-tiles by their scale relative to $H$ and by the translate of $H$ containing the center of their spatial interval: for $\nu\geq1$ and $k\in\Z$, let $H_k\coloneqq H+k\ell(H)$ and
\[
\begin{split}
 &\bS^{\nu,k}_H\coloneqq\big\{s\in\bS:\ \ell(I_s)=2^{-\nu}\ell(H),\ c(I_s)\in H_k\big\},\\
 &\Psi^{\nu,k}_H\coloneqq\sum_{s\in\bS^{\nu,k}_H}\overline{\beta_s}\,\varphi_{s_1},\qquad n^{\nu,k}_H\coloneqq\#\bS^{\nu,k}_H .
\end{split}
\]
As the lengths of dyadic intervals are powers of $2$, and the translates $H_k$, $k\in\Z$, partition $\R$, the collection $\bS^{\mathsf{fine}}(H)$ is the disjoint union of the collections $\bS^{\nu,k}_H$, so that
\[
 \Psi^{\mathsf{fine}}_H=\sum_{\nu\geq1}\sum_{k\in\Z}\Psi^{\nu,k}_H .
\]
By \eqref{e:finesep}, $\bS^{\nu,k}_H=\varnothing$ if $|k|\leq3$. If $|k|\geq4$, the distance of $3H$ from $H_k$ is at least $\frac{|k|}2\ell(H)$, that is, $2^{\nu-1}|k|$ times the length of the spatial intervals of $\bS^{\nu,k}_H$. Hence, by \eqref{e:chi}, the wave packets of $\Psi^{\nu,k}_H$ have relative size $(2^\nu|k|)^{-M}$ on $3H$, 
\begin{equation}\label{e:finetail}
 |\varphi_{s_1}(x)|\lesssim\frac{(2^\nu|k|)^{-M}}{|I_s|},\qquad x\in3H,\ s\in\bS^{\nu,k}_H .
\end{equation}
The next lemma provides the analogues of \eqref{e:FH} and \eqref{e:count} for the functions $\Psi^{\nu,k}_H$. The losses, polynomial in $2^\nu|k|$, are compensated by the decay in \eqref{e:finetail}.

\begin{lemma}\label{l:fine}
Let $\nu\geq1$ and $|k|\geq4$. Then
\begin{align}
 \label{e:fine1}\|\Psi^{\nu,k}_H\|_{L^\infty(3H)}&\lesssim(2^\nu|k|)^{-M}\,n^{\nu,k}_H,\qquad H\in\cW,\\
 \label{e:fine2}\sum_{H\in\cW}|H|\,n^{\nu,k}_H&\lesssim2^{2\nu}|k|\,\cN,\\
 \label{e:fine3}\sum_{H\in\cW}\|\Psi^{\nu,k}_H\|_{L^2(3H)}^2&\lesssim2^{2\nu}(2^\nu|k|)^{-2M}\,\cN .
\end{align}
\end{lemma}
\begin{proof}
The bound \eqref{e:fine1} follows from \eqref{e:finetail} and $|\beta_s|\leq|I_s|$.

To prove \eqref{e:fine2}, fix $\bT\in\cF$, and let $H\in\cW$ with $\bT\cap\bS^{\nu,k}_H\neq\varnothing$. The spatial intervals of the tri-tiles of $\bT\cap\bS^{\nu,k}_H$ have length $2^{-\nu}\ell(H)$ and centers in $H_k$, so that there are at most $2^\nu$ of them, each shared by at most $2^\kappa$ tri-tiles of $\bT$. Moreover, $\ell(H)\leq2^\nu\ell_\bT$ and $H_k$ intersects $I_\bT$, so that $H$ is contained in the set of points at distance at most $(|k|+1)2^\nu\ell_\bT$ from $I_\bT$. As the intervals of $\cW$ are pairwise disjoint,
\[
 \sum_{H\in\cW}|H|\,\#\big(\bT\cap\bS^{\nu,k}_H\big)\leq2^{\kappa+\nu}\sum_{\substack{H\in\cW\\ \bT\cap\bS^{\nu,k}_H\neq\varnothing}}|H|\lesssim2^{2\nu}|k|\,|I_\bT| .
\]
Summing over $\bT\in\cF$ proves \eqref{e:fine2}.

To prove \eqref{e:fine3}, fix $H$, let $\ell\coloneqq2^{-\nu}\ell(H)$, and split $\Psi^{\nu,k}_H=\sum_\omega\Psi^{\nu,k}_{H,\omega}$ according to the value $\omega$ of $\omega_{s_1}$. Fix a nonzero Schwartz function $\rho\geq0$ with $\supp\widehat\rho\subset[-1,1]$, and let
\[
 \zeta\coloneqq c\,\ind_{4H}*\rho_\ell,\qquad\rho_\ell\coloneqq\ell^{-1}\rho(\cdot/\ell),
\]
with $c>0$ depending only on $\rho$ such that $\zeta\geq1$ on $3H$, which is possible as $\ell\leq\ell(H)$. Then
\[
 \supp\widehat\zeta\subset[-\ell^{-1},\ell^{-1}],\qquad\|\zeta\|_2^2\lesssim|H|,\qquad|\zeta(x)|\lesssim\Big(1+\frac{\dist(x,4H)}\ell\Big)^{-M}.
\]
The functions $\zeta\Psi^{\nu,k}_{H,\omega}$ have frequency supports in the intervals $\omega+[-\ell^{-1},\ell^{-1}]$, which have bounded overlap, so that
\[
 \|\Psi^{\nu,k}_H\|_{L^2(3H)}^2\leq\|\zeta\Psi^{\nu,k}_H\|_2^2\lesssim\sum_\omega\|\zeta\Psi^{\nu,k}_{H,\omega}\|_2^2 .
\]
Let $s\in\bS^{\nu,k}_H$. As $\dist(4H,H_k)\gtrsim|k|\ell(H)$, every $x\in\R$ satisfies $\dist(x,4H)\gtrsim2^\nu|k|\ell$ or $|x-c(I_s)|\gtrsim2^\nu|k|\ell$. Splitting the integral accordingly, \eqref{e:chi} and $\|\zeta\|_2^2\lesssim|H|=2^\nu|I_s|$ give
\[
 \|\zeta\varphi_{s_1}\|_2^2\lesssim\frac{2^\nu(2^\nu|k|)^{-2M}}{|I_s|} .
\]
For fixed $\omega$, the intervals $I_s$ with $\omega_{s_1}=\omega$ are distinct, by the injectivity of $s\mapsto s_1$; as they have length $\ell$ and centers in $H_k$, there are at most $2^\nu$ of them. Hence, by the Cauchy--Schwarz inequality,
\[
 \|\zeta\Psi^{\nu,k}_{H,\omega}\|_2^2\leq2^\nu\sum_{\omega_{s_1}=\omega}|\beta_s|^2\,\|\zeta\varphi_{s_1}\|_2^2\lesssim2^{2\nu}(2^\nu|k|)^{-2M}\sum_{\omega_{s_1}=\omega}\frac{|\beta_s|^2}{|I_s|} .
\]
Finally, for fixed $s$, at most one $H\in\cW$ satisfies $s\in\bS^{\nu,k}_H$, namely the interval of $\cW$ of length $2^\nu\ell(I_s)$ containing $c(I_s)-k2^\nu\ell(I_s)$. Summing over $\omega$ and $H$, \eqref{e:fine3} follows from \eqref{e:coef}.
\end{proof}

We now conclude the proof of \eqref{e:badfine}. As $\bS^{\nu,k}_H=\varnothing$ for $|k|\leq3$ and $b_H=f\ind_H-g_H$,
\[
 \langle b_H,\Psi^{\mathsf{fine}}_H\rangle=\sum_{\nu\geq1}\sum_{|k|\geq4}\Big(\langle f\ind_H,\Psi^{\nu,k}_H\rangle-\langle g_H,\Psi^{\nu,k}_H\rangle\Big),\qquad H\in\cW .
\]
Hence it suffices to bound, for $\nu\geq1$ and $|k|\geq4$, the sums over $H\in\cW$ of the absolute values of the two terms in parentheses by a constant times $2^{-\nu}|k|^{-2}\cN^{\frac{1+\alpha}2}$. For the first term, we apply Lemma~\ref{l:sum} on $O=E_1$, as in Section~\ref{ss:coarse}, to $u=f$ and to
\[
 F\coloneqq\sum_{H\in\cW}\Psi^{\nu,k}_H\ind_H,\qquad\Upsilon\coloneqq(2^\nu|k|)^{-M}\sum_{H\in\cW}n^{\nu,k}_H\ind_H .
\]
By \eqref{e:fine1}, $|F|\lesssim\Upsilon$ on $E_1$. By \eqref{e:fine2} and \eqref{e:fine3},
\[
 \|\Upsilon\|_{L^1(E_1)}\lesssim2^{-\nu(M-2)}|k|^{1-M}\,\cN,\qquad\|F\|_{L^2(E_1)}^2\lesssim2^{-2\nu(M-1)}|k|^{-2M}\,\cN .
\]
As $|E_1|\lesssim1$ and $M$ is large, Lemma~\ref{l:sum} and \eqref{e:whitney} give
\[
 \sum_{H\in\cW}|\langle f\ind_H,\Psi^{\nu,k}_H\rangle|\leq\|fF\|_{L^1(E_1)}\lesssim2^{-\nu(M-2)}|k|^{1-M}\,\cN^{\frac{1+\alpha}2}\leq2^{-\nu}|k|^{-2}\,\cN^{\frac{1+\alpha}2} .
\]
For the second term, the Cauchy--Schwarz inequality, \eqref{e:good} and \eqref{e:fine3} give
\[
\begin{split}
 \sum_{H\in\cW}|\langle g_H,\Psi^{\nu,k}_H\rangle|&\leq\Big(\sum_{H\in\cW}\|g_H\|_{L^2(3H)}^2\Big)^{\frac12}\Big(\sum_{H\in\cW}\|\Psi^{\nu,k}_H\|_{L^2(3H)}^2\Big)^{\frac12}\\
 &\lesssim2^{-\nu(M-1)}|k|^{-M}\,\cN^{\frac{1+\alpha}2}\leq2^{-\nu}|k|^{-2}\,\cN^{\frac{1+\alpha}2} .
\end{split}
\]
This proves \eqref{e:badfine}, and completes the proof of Lemma~\ref{l:bad}.

Finally, inserting \eqref{e:goodterms} and Lemma~\ref{l:bad} into \eqref{e:split} yields
\[
 |\Lambda_\bS(f,f_2,f_3)|\lesssim\max\{1,\cN\}^{\frac{1+\alpha}2},
\]
which is Proposition~\ref{p:level} under the normalization \eqref{e:norm37}. \qed

\section{Remarks}\label{s:remarks}

\subsection{Uniformity} The coefficients $\eps_s$ and the wave packets $\varphi_{s_j}$ enter the proofs only through $|\eps_s|\leq1$ and $\varphi_{s_j}\in\Phi(s_j)$, and the major subset $F_3'$ of Proposition~\ref{p:main} depends neither on $m$ nor on the model sums. Hence Proposition~\ref{p:main} and Theorem~\ref{t:main} hold uniformly over the class, which is the form in which Theorem~\ref{t:main} is sharp.

\subsection{The order of the moments at the tops} With moments of order zero only, as in \cites{NOT,DPT}, the remainder $\Delta_s$ in the proof of Lemma~\ref{l:near} is of relative size $2^{-\mu}$ instead of $2^{-2\mu}$, so each relative scale contributes a bounded amount to each bin, and the bound of Lemma~\ref{l:near} becomes the number of relative scales, which is unbounded. This is the loss which \cite{DPT} compensates with additional frequencies near each top and with wave packets of almost exponential decay. Moments of order one at the tops are the smallest choice which makes the sum over the relative scales geometric; no higher order and no separation of scales is needed for the near trees.

\subsection{Separation of scales} In \cite{DPOu}, where an embedding of the wave packet transform is proved, the bad part must have small size on all tents, including those at the scale of its support; this is achieved with moments of high order and a dilation of the exceptional set by a small power of the oscillation of the counting function, which enters the final estimate as a loss in the exponent. Here only the sum of the interactions over each Whitney interval is estimated, through the orthogonality of the bins and the pairing with $f_1\ind_H$, so the tri-tiles of scale comparable to $\ell(H)$ need not be excised, and no dilation of the exceptional set occurs.

\subsection{Far trees and the counting function} The dimension spent on the far trees is controlled by $n^{\mathsf{far}}_H$, a weighted count of the far tops in which each tree carries a weight decaying with its relative distance from $H$ and independent of the number of its scales; by \eqref{e:count}, its integral over $E_1$ is controlled by $\|\sum_{\bT}\ind_{I_\bT}\|_1$. In particular no dilation of the tops by a growing factor and no bound on the oscillation of the counting function is used. Merging of near tops in a bin costs nothing, by the counting of Lemma~\ref{l:near}; merging of far tops in a bin costs $O(1+\log n^{\mathsf{far}}_{H,B})$ dimensions, by Lemma~\ref{l:far}.

\subsection*{Declaration of generative AI and AI-assisted technologies
in the manuscript preparation process}
The author used Claude
(Anthropic) to assist in streamlining a few  proofs from the author's original draft, in
checking the mathematical arguments, cross-references and citations,
and in preparing and editing the \LaTeX{} source. The author reviewed
and edited the manuscript, independently verified the mathematical
arguments and references, and takes full responsibility for the
content of the article.

\bibliographystyle{plain}
\bibliography{refs}
\end{document}